\magnification = \magstephalf

\newcount\footnum
\footnum=1
\def\newftnum{\global\advance\footnum by 1}
\def\ft#1{\footnote{$^{\vrule width .2cm\the\footnum\vrule width .2cm}$}{#1}
\newftnum}

\newcount\fnum
\fnum=1
\def\newfnum{\global\advance\fnum by 1}
\def\ftt#1{\footnote{$^{\clubsuit\the\fnum \clubsuit}$}{#1}\newfnum}



\newcount\sectno
\newcount\subsectno
\newcount\parno
\newcount\equationno
\newif\ifsubsections
\subsectionsfalse

\def\sectnum{\the\sectno}
\def\subsectnum{\sectnum\ifsubsections .\the\subsectno\fi}
\def\parnum{\subsectnum .{\the\parno}}
\def\eqnum{\subsectnum .\the\equationno}

\def\Title#1{
        \sectno=0
\centerline{\ifproofmode\rmarginsert[\timestamp]\fi\titlefont #1}
        \vskip .75pc}

\def\abstract#1\endabstract
{
{\abstractfont
       \baselineskip=9pt
       \leftskip=4pc  \rightskip=4pc
       \bigskip
       \noindent
        ABSTRACT.\ #1
\medskip}
}

\def\thanks[#1]#2\endthanks{\footnote{$^#1$}{\footnotefont\kern-6pt #2}}

\newcount\minutes
\newcount\scratch

\def\timestamp{%
\scratch=\time
\divide\scratch by 60
\edef\hours{\the\scratch}
\multiply\scratch by 60
\minutes=\time
\advance\minutes by -\scratch
\the \month/\the\day$\,$---$\,$\hours:\null
\ifnum\minutes< 10 0\fi
\the\minutes}

\def\today{\ifcase\month\or
January\or February\or March\or April\or May\or June\or
July\or August\or September\or October\or November\or December\fi
\space\number\day,\space\number\year}

\outer\def\newsection #1.\par{\vskip1.5pc plus.75pc \penalty-250
        \subsectno=0
        \parno=0
        \equationno=0
        \advance\sectno by1
        \leftline{\smalltitlefont \sectnum.\hskip 1pc  #1}
                   \nobreak \vskip.75pc\noindent}

\outer\def\newsectiontwoline #1/#2/.\par{\vskip1.5pc plus.75pc \penalty-250
        \subsectno=0
        \parno=0
        \equationno=0
        \advance\sectno by1
        \leftline{\smalltitlefont \sectnum.\hskip 1pc  #1}
        \leftline{\smalltitlefont \hskip 22pt #2}
                   \nobreak \vskip.75pc\noindent}

\outer\def\newsubsection #1.\par{\vskip1pc plus.5pc\penalty-250
        \parno=0
        \equationno=0
        \advance\subsectno by1
        \leftline{{\bf \subsectnum}\hskip 1pc  #1.}
                   \nobreak \vskip.5pc\noindent}

\def\newpar #1.{\advance \parno by1
        \par
    \medbreak \noindent
         {\bf \parnum. #1.} \hskip 6pt}

\long\def \newclaim #1. #2\par {\advance \parno by1
       \medbreak \noindent
        {\bf \parnum \hskip 6 pt #1.\hskip 6pt} {\sl #2} \par \medbreak}

\def\eq $$#1$${\global \advance \equationno by1 $$#1\eqno(\eqnum)$$}

\def\rmarginsert[#1]{\hglue 0pt\vadjust
{\null\vskip -\baselineskip\rightline{\abstractfont\rlap{\hfil\  #1}}}}

\def\lmarginsert[#1]{\hglue 0pt\vadjust
{\null\vskip -\baselineskip\leftline{\abstractfont\llap{#1\ \hfill}}}}

\newif\ifproofmode
\proofmodefalse

\def\refpar[#1]#2.{\advance \parno by1
        \par
    \medbreak \noindent
         {\bf \parnum \hskip 6 pt #2.\hskip 6pt}%
\expandafter\edef\csname ref#1\endcsname
{\parnum}\ifproofmode\rmarginsert[\string\ref#1]\fi}

\long\def \refclaim[#1]#2. #3\par {\advance \parno by1
       \medbreak \noindent
{\bf \parnum \hskip 6 pt #2.\hskip 6pt}%
\expandafter\edef\csname ref#1\endcsname
{\parnum}\ifproofmode\rmarginsert[\string\ref#1]\fi
{\sl #3} \par \medbreak}

\def\refer[#1]{%
\expandafter\xdef\csname ref#1\endcsname
{\parnum}\ifproofmode\rmarginsert[\string\ref#1]\fi}

\def\refereq[#1]$$#2$$ {%
\eq$$#2$$%
\expandafter\xdef\csname ref#1\endcsname{(\eqnum)}%
\ifproofmode\rmarginsert[\string\ref#1]\fi
}

\def\refeq{\refereq}

\def \Definition #1\\ {\vskip 1pc \noindent
         {\bf #1. Definition. \hskip 6pt}\vskip 1pc}

\def\proof{{PROOF.} \enspace}

\def\qedmark{\hbox{\vrule height 4pt width 3pt}}
\def\qedskip{\vrule height 4pt width 0pt depth 1pc}
\def\qed{\penalty 1000\quad\penalty 1000{\qedmark\qedskip}}

\def \a {\alpha}
\def \b {\beta}
\let\C=\Cx
\def \d {\delta}

\def \g {\gamma}

\def \K {\nabla}
\def \l {\lambda}
\def \L {\Lambda}
\def \n {\,\vert\,}
\def \N {\,\Vert\,}
\def \o {\theta}
\mathchardef\p="011E    
\def \R{\reals}
\def \s {\sigma}

\def\II{I\!I}     

\def \Gr {{\mbi Gr}}

\def\Gtwo{{\mathop{{{\mbi G\/}}\kern-.5pt_{{}_2}}}}
\def\Ffour{{\mathop{{{\mbi F\/}}\kern-2.5pt_{{}_4}}}}
\def\Esix{{\mathop{{{\mbi E\/}}\kern-.5pt_{{}_6}}}}
\def\Eseven{{\mathop{{{\mbi E\/}}\kern-.5pt_{{}_7}}}}
\def\Eeight{{\mathop{{{\mbi E\/}}\kern-.5pt_{{}_8}}}}


\def\cc{{\liefont C}}

\def\cg{{\liefont G}}

\def\ck{{\liefont K}}

\def\cm{{\liefont M}}

\def\co{{\liefont O}}

\def\cp{{\liefont P}}

\def\cs{{\liefont S}}

\def\ct{{\liefont T}}

\def\cu{{\liefont U}}

\def\sdp{
\mathop{\hbox{$\raise 1pt\hbox{$\scriptscriptstyle |$}\kern-2.5pt\times$}}
           }

\def\dag{{\raise 1 pt\hbox{{$\scriptscriptstyle \dagger$}}}}

\def\*{\raise 1.5pt \hbox{*}}

\def\tr{\mathop{\rm tr}\nolimits}

\def\diag{\mathop{\rm diag}\nolimits}

\def\Ad{\mathop{\rm Ad}\nolimits}
\def\ad{\mathop{\rm ad}\nolimits}

\def\diag {\mathop{\rm diag}\nolimits}

\def\id {\mathop{\rm id}\nolimits}

\def \li{\langle}
\def \ri{\rangle}


\def\ifundefined#1{\expandafter\ifx\csname#1\endcsname\relax}

\def\cross{\times}
\def\longerrightarrow{-\kern-5pt\longrightarrow}

\def\star{\lower 1pt\hbox{*}}
\def \nulset {
\raise 1pt\hbox{
\hskip -3pt$\not$\kern -0.2pt \raise .7pt\hbox{${\scriptstyle\bigcirc}$}}}
\def \norm|#1|{\Vert#1\Vert}

\def \interior(#1){#1\kern -6pt \raise 7.5pt
         \hbox{$\scriptstyle \circ$}{}\hskip 2pt}

\def\twist_#1{\kern -.15em\cross\kern -.30 em{}_{{}_{#1}}\kern .07 em}

\font\cmr=cmr10 at 10pt
\font\cmrviii=cmr8
\font\cmrvi=cmr6

\font\cmrXIV=cmr12 at 14 pt
\font\cmrXX=cmr17 at 20 pt
\font\cmrXXIV=cmr17 at 24 pt
\font\cmbxXII=cmbx12
\font\cmbxsl=cmbxsl10
\font\cmbxslviii=cmbxsl10 at 8pt
\font\cmbxslv=cmbxsl10 at 5pt

         \font\tenrm=cmr10 at 10.3 pt
         \font\sevenrm=cmr7 at 7.21 pt
         \font\fiverm=cmr5 at 5.15 pt
         \font\teni=cmmi10 at 10.3 pt
         \font\seveni=cmmi7 at 7.21 pt
         \font\fivei=cmmi5 at 5.15 pt
         \font\tensy=cmsy10 at 10.3 pt
         \font\sevensy=cmsy7 at 7.21 pt
         \font\fivesy=cmsy5 at 5.15 pt
         \font\tenex=cmex10 at 10.3 pt
         \font\tenbf=cmbx10 at 10.3 pt
         \font\sevenbf=cmbx7 at 7.21 pt
         \font\fivebf=cmbx5 at 5.15 pt

\def\UseComputerModern   
{
\textfont0=\tenrm \scriptfont0=\sevenrm \scriptscriptfont0=\fiverm
\def\rm{\fam0\tenrm}
\textfont1=\teni \scriptfont1=\seveni \scriptscriptfont1=\fivei
\def\mit{\fam1} \def\oldstyle{\fam1\teni}
\textfont2=\tensy \scriptfont2=\sevensy \scriptscriptfont2=\fivesy
\def\cal{\fam2}
\textfont3=\tenex \scriptfont3=\tenex \scriptscriptfont3=\tenex
\def\it{\fam\itfam\tenit} 
\textfont\itfam=\tenit
\def\sl{\fam\slfam\tensl} 
\textfont\slfam=\tensl
\def\bf{\fam\bffam\tenbf} 
\textfont\bffam=\tenbf \scriptfont\bffam=\sevenbf
\scriptscriptfont\bffam=\fivebf
\def\tt{\fam\ttfam\tentt} 
\textfont\ttfam=\tentt
\def\abstractfont{\cmrviii}
\def\footnotefont{\cmrviii}
\def\tinyfont{\cmrvi}
\def\smalltitlefont{\cmbxXII}  
\def\titlefont{\cmrXIV}
\def\bigtitlefont{\cmrXX}
\def\verybigtitlefont{\cmrXXIV}
\textfont9=\cmbxsl \scriptfont9=\cmbxslviii \scriptscriptfont9=\cmbxslv
\def\mbi{\fam9}
\cmr
}

\def\liefont{\cal}
\font\bbb=msbm10

\def\R{\hbox {\bbb R} }
\def\C{\hbox {\bbb C} }

\def \bs {\bigskip}
\def \ms {\medskip}
\def \ss {\smallskip}

\def \ni {\noindent}

\def\enditem{\item{}\par\vskip-\baselineskip\noindent}
\def\ei{\enditem}

\def\id{\mathop{\rm id}\nolimits}

\def\ni{\noindent}
\def\bs{\bigskip}
\def\ms{\medskip}

\def\ss{\smallskip}
\def\ct{{\cal T}}

      \baselineskip= 11  pt
      \hsize  6.5 true in 
      \vsize  9.1true  in

\UseComputerModern
\subsectionsfalse
\font\smalltitlefont=cmbx10 at 11 pt
%


\def\Bibliography
{

\font\TRten=cmr10 at 10 true pt
\font\TIten=cmti10 at 10 true pt
\font\TBten=cmbx10 at 10 true pt

\def\ourindent{\hfil\vskip-\baselineskip}

    \frenchspacing
    \parindent=0pt

    \def \keyfnt{\TRten}
    \def \authornamefnt{\TRten}
    \def \booktitlefnt{\TIten}
    \def \articletitlefnt{\TRten}
    \def \journalnamefnt{\TIten}
    \def \volumefnt{\TBten}
    \def \publishernamefnt{\TRten}
    \def \pagesfnt{\TRten}
    \def \yearfnt{\TRten}
    \def \commentfnt{\TRten}

    \def \bookitem //##1//##2//##3//##4//##5//##6//##7//##8//
         { \goodbreak{\par\hskip-40pt{\keyfnt
[##1]}\ourindent{\authornamefnt ##2,}}
                 {\booktitlefnt ##3.\/}\thinspace
                 {\publishernamefnt ##4,}
                 {\yearfnt ##6.}
                 {\commentfnt ##8}
          }

\def \b{\bookitem}

\def \articleitem //##1//##2//##3//##4//##5//##6//##7//##8//

         { \goodbreak{\par\hskip-40pt{\keyfnt
[##1]}\ourindent{\authornamefnt ##2,}}
                 {\articletitlefnt ##3},
                 {\journalnamefnt ##4\/}
                 {\volumefnt ##5}
                 {\hbox{\yearfnt(\hskip -1pt ##6)}},
                 {\pagesfnt ##7.}
                 {\commentfnt ##8}
          }
\def \a{\articleitem}

\def \preprintitem //##1//##2//##3//##4//##5//##6//##7//##8//
         { \goodbreak{\par\hskip-40pt{\keyfnt
[##1]}\ourindent{\authornamefnt ##2,}}
                 {\articletitlefnt ##3},
                 {\commentfnt ##8}
          }
\def \p{\preprintitem}

      \vskip 1in
      \centerline{References}
      \vskip .5in
}


\def\Reduce#1
{
\font\stenrm=cmr10 scaled #1
\font\sninerm=cmr9 scaled #1
\font\seightrm=cmr8 scaled #1
\font\ssevenrm=cmr7 scaled #1
\font\ssixrm=cmr6 scaled #1
\font\sfiverm=cmr5 scaled #1

\font\steni=cmmi10 scaled #1
\font\sninei=cmmi9 scaled #1
\font\seighti=cmmi8 scaled #1
\font\sseveni=cmmi7 scaled #1
\font\ssixi=cmmi6 scaled #1
\font\sfivei=cmmi5 scaled #1

\font\stenit=cmti10 scaled #1
\font\snineit=cmti9 scaled #1
\font\seightit=cmti8 scaled #1
\font\ssevenit=cmti7 scaled #1

\font\stensy=cmsy10 scaled #1
\font\sninesy=cmsy9 scaled #1
\font\seightsy=cmsy8 scaled #1
\font\ssevensy=cmsy7 scaled #1
\font\ssixsy=cmsy6 scaled #1
\font\sfivesy=cmsy5 scaled #1

\font\stenbf=cmbx10 scaled #1
\font\sninebf=cmbx9 scaled #1
\font\seightbf=cmbx8 scaled #1
\font\ssevenbf=cmbx7 scaled #1
\font\ssixbf=cmbx6 scaled #1
\font\sfivebf=cmbx5 scaled #1

\font\stentt=cmtt10  scaled #1
\font\sninett=cmtt9 scaled #1
\font\seighttt=cmtt8 scaled #1

\font\stenex=cmex10 scaled #1

\font\stensl=cmsl10 scaled #1
\font\sninesl=cmsl9 scaled #1
\font\seightsl=cmsl8 scaled #1

\textfont0=\stenrm \scriptfont0=\ssevenrm \scriptscriptfont0=\sfiverm
\def\rm{\fam0\stenrm}
\textfont1=\steni \scriptfont1=\sseveni \scriptscriptfont1=\sfivei
\def\mit{\fam1} \def\oldstyle{\fam1\steni}
\textfont2=\stensy \scriptfont2=\ssevensy \scriptscriptfont2=\sfivesy
\def\cal{\fam2}
\textfont3=\stenex \scriptfont3=\stenex \scriptscriptfont3=\stenex
\def\it{\fam\itfam\tenit} 
\textfont\itfam=\stenit
\def\sl{\fam\slfam\stensl} 
\textfont\slfam=\stensl
\def\bf{\fam\bffam\stenbf} 
\textfont\bffam=\stenbf \scriptfont\bffam=\ssevenbf
\scriptscriptfont\bffam=\sfivebf
\def\tt{\fam\ttfam\stentt} 
\textfont\ttfam=\stentt
\rm
}

\def\p{\partial}

\def\p{\partial}

\def\II{\hbox{\rm II\/}}
\def\Id{\hbox{\rm Id\/}}
\def\ti{\tilde}

      \hsize  6.5 true in 
\voffset=  1 true in
      \vsize   7.875 true  in

\Title{Completely integrable curve flows on Adjoint orbits}

\bs
\centerline{Chuu-Lian Terng\footnote{$^1$}{Research supported in part by
NSF Grant DMS 9972172 and Humboldt Senior Scientist Award.} and Gudlaugur
Thorbergsson\footnote{$^2$}{Research supported in part by the Deutsche
Forschungsgemeinschaft.}}
\bs\bs
\centerline{Dedicated to Professor S. S. Chern on his 90th Birthday\/}
\bs
\centerline{ Abstract}
\bs

It is known that the Schr\"odinger flow on a complex
Grassmann manifold is equivalent to the matrix non-linear Schr\"odinger equation and the
Ferapontov flow on a principal Adjoint $U(n)$-orbit is equivalent to the
$n$-wave equation.  In this paper, we give a systematic
method to construct integrable geometric curve flows on Adjoint $U$-orbits from flows in the
soliton hierarchy associated to a compact Lie group $U$. There are natural geometric
bi-Hamiltonian structures on the space of curves on Adjoint orbits, and they correspond to the
order two and three Hamiltonian structures on soliton equations under our
construction.  We study the Hamiltonian theory of these geometric curve flows and also give
several explicit examples.

\bs
\ni {\it Keywords\/}: completely integrable Hamiltonian systems, curve flows, symmetric
spaces. 

\ni {\it 2000 Mathematics Subject Classification number\/}: 37K10, 37K25, 53C35
\vfil\eject

      \hsize  6.5 true  in 
\voffset=- 10 true pt
      \vsize   9.1 true  in

\newsection  Introduction.\par

  There
are several natural geometric flows on Adjoint
orbits
that are  known to be equivalent to soliton equations. The first example
is the Heisenberg ferromagnetic model (HFM) for
$\g:\R^2\to S^2$,
$$\g_t=\g\times \g_{xx}, \eqno({\rm HFM\/})$$ where
$\times$ is the cross product in $\R^3$.
    It was proved in [FT] that the HFM
is equivalent to  the non-linear Schr\"odinger equation
(NLS):
    $$q_t= (q_{xx}+2 \n q\n^2 q). \eqno({\rm NLS\/})$$

The second example is the Schr\"odinger flow on the Hermitian symmetric space
$\Gr(k,\C^n)$.  Recall that the Schr\"odinger flow on a K\"ahler
manifold $M$ is the
evolution equation on the space of maps from $\R$ to $M$:
$$\g_t= j_\g(\K_{\g_x}\g_x),$$
where $j$ is the complex structure and $\K$ is the Levi-Civita
connection of the
K\"ahler metric on $M$.  The
Adjoint
$U(n)$-orbit
$M$  at
\refeq[hq]$$a={1\over 2} \pmatrix{i\ \Id_k&0\cr 0 & -i\
\Id_{n-k}\cr}$$ equipped with the
induced metric from the inner product $\li u_1,u_2\ri = -\tr(u_1u_2)$ on $u(n)$
    is isometric to the Hermitian symmetric space $\Gr(k,\C^n)$, and
$j_x=\ad(x)$ is the complex structure on $\Gr(k,\C^n)$.
     Terng and Uhlenbeck showed in [TU2] that the Schr\"odinger
flow on $\Gr(k,\C^n)$ is
\refeq[hr]$$\g_t=j_\g(\K_{\g_x}\ \g_x) =[\g, \g_{xx}],$$
and is equivalent to the matrix non-linear Schr\"odinger
equation (MNLS) for  maps $q$ from
$\R^2$ to the space
$\cm_{k\times (n-k)}$ of $k\times (n-k)$ complex matrices:
$$q_t= (q_{xx} + 2 qq^* q),\eqno({\rm MNLS\/})$$
where $q^*=\bar q^t$.
Note that when $n=2$, $M$ is isometric to the round sphere. If we
identify $su(2)$ with
$\R^3$ in the usual way, then $[\xi,\eta]$ corresponds to the cross
product $\xi\times \eta$. So
the Schr\"odinger flow
\refhr{} on
$S^2$ is the HFM.

The third example is due to
Ferapontov. Let  $\cu$ denote the Lie algebra
of a compact Lie group $U$.   An element $a\in \cu$ is called {\it
regular\/} ({\it singular\/} respectively) if the Adjoint
$U$-orbit $M_a$ at $a$ in
$\cu$ is a principal (singular respectively) orbit.
     Let $\li \ , \ \ri$ denote an Ad-invariant inner product on $\cu$,
$\cu_a$ the isotropy
subalgebra of $a$, and  $\cu_a^\perp$ the orthogonal complement of
$\cu_a$ in $\cu$
with respect to $\li\ , \ \ri$.   If $a\in \cu$ is regular, then it
is known that ([Te1])
\item {$\diamond$} the normal
bundle $\nu(M_a)$ is flat,
\item {$\diamond$} $\nu(M_a)_x=\cu_x$, which is the maximal
abelian subalgebra of $\cu$ containing $x$,
\item {$\diamond$} given any $b\in \cu_a$, the map $\hat b$
defined by
\refeq[eo]$$\hat b(gag^{-1})= gbg^{-1}, \quad g\in U,$$
is a well-defined parallel normal field of $M_a$.

\ni Ferapontov proved in [F4] that a solution $u:\R^2\to \cu_a^\perp$ 
of the {\it
n-wave equation associated to $\cu$\/},
\refeq[ab]$$u_t= \ad(b)\ad(a)^{-1}(u_x) + [u, \ad(b)\ad(a)^{-1}(u)],$$
gives rise to a solution of the following curve flow on $M_a$:
\refeq[el]$$\g_t= (\hat b(\g))_x = -A_{\hat b(\g)}(\g_x).$$
Here $A_v$ is the shape operator along $v$.

The MNLS and the $n$-wave equation
\refab{} are flows in the
$U(n)$- and
$U$-hierarchy of soliton flows respectively.
The $U$-hierarchy of soliton flows is obtained by restricting the
ANKS-ZS hierarchy to an invariant submanifold associated to the
reality condition given by
the real group $U$.  One goal of this paper is to give a systematic method to
construct geometric curve flows on an
Adjoint
$U$-orbit for each flow in  the $U$-hierarchy
that includes all three examples given above.

We review the
construction of the $U$-hierarchy
next.   Let
$a\in \cu$, and $\cs(\R,\cu_a^\perp)$ the space of maps from
$\R$ to $\cu_a^\perp$ that decay rapidly at infinity.
The $U$-hierarchy defined by $a$ is a collection of commuting Hamiltonian flows
on $\cs(\R,\cu_a^\perp)$ (cf.\ [TU1]).   When $a$ is
regular, the
$U$-hierarchy is parametrized by $(b,j)$ with
$b\in \cu_a$ and  $j$ positive integer.  The $(b,j)$-flow  is
\refeq[co]$$u_t= (Q_{b,j}(u))_x+ [u,Q_{b,j}(u)],$$
where $Q_{b,j}(u)$ are $\cu$-valued maps  determined by
the following conditions:
\refeq[ad]$$(Q_{b,j}(u))_x + [u, Q_{b, j}(u)] = [Q_{b,j+1}(u),a],
\quad Q_{b,0}(u)=b,$$
and
\refeq[iw]$$\sum_{j=0}^\infty Q_{b,j}(u)\l^{-j} \ {\rm is\ conjugate
\ to \ } b \
{\rm as\ an \ asymptotic\ expansion\/}.$$
These conditions imply that $Q_{b,j}(u)$ is a polynomial in
$u, \p_x u, \ldots, \p_x^{j-1}u$ and for $j\geq 1$
\refeq[jd]$$ Q_{b,j}(u)\in \cs(\R, \cu) \quad {\rm if \ } u\in
\cs(\R,\cu_a^\perp).$$
For more detail see [Sa] and [TU1]. 

  When $a$ is singular, the
$U$-hierarchy is the
collection of
$(a,j)$-flows:
$$u_t=(Q_{a,j}(u))_x + [u,Q_{a,j}(u)].$$

Recall that a $\cg$-valued connection 1-form $w=Adx + Bdt$ is {\it flat\/} if
$dw=-w\wedge w$ or
equivalently,
$$A_t-B_x = [A,B].$$
     The recursive formula \refad{} implies that $u$ is a solution of the
$(b,j)$-flow \refco{} if and only if
\refeq[ex]$$\o_\l=(a\l + u)\ dx + (b\l^j + Q_{b,1}(u)\l^{j-1} +
\cdots + Q_{b,j}(u))\ dt$$
is a flat $\cu_{\C}$-valued connection 1-form over the $(x,t)$ plane
for all $\l\in
\C$.  The $1$-form $\o_\l$ is called a {\it Lax pair\/} of the
$(b,j)$-flow.

\ms
For example:

\item {$\diamond$} If $U=SU(2)$ and
$a={1\over 2}\diag(i,-i)$, then
$\cu_a^\perp=\left\{\pmatrix{0&q\cr -\bar q &0\cr}\bigg| q\in \C\right\}$.
Identify $\cs(\R, \cu_a^\perp)$ with $\cs(\R, \C)$. The $(a,2)$-flow in the
$SU(2)$-hierarchy  is the NLS.

\item {$\diamond$} Let $a\in U$ be a regular element,
and $b\in \cu_a$.  The $(b,1)$-flow in
the
$U$-hierarchy on
$\cs(\R,
\cu_a^\perp)$ is the
$n$-wave equation \refab{}.

\item {$\diamond$} Let $a\in u(n)$ be as in \refhq{}. Then
$\cu_a^\perp$ is the space of $\pmatrix{0& q\cr -q^*& 0\cr}$ with $q$ a
$k\times (n-k)$ complex matrix.   Identify $u$ with
$q$.  The
$(a,2)$-flow of the
$U(n)$-hierarchy on $\cs(\R,\cu_a^\perp)$ is the MNLS.

\item {$\diamond$}  If $U/K$ is a Hermitian symmetric space, then
(cf.~[H], [W]) there exists
$a\in \ck$ such that $\cu_a=\ck$, $\ad(a)^2=-\Id$ on $\cu_a^\perp$, and
$\cu=\cu_a+\cu_a^\perp$ is a Cartan decomposition. The Adjoint
$U$-orbit at $a$ in
$\cu$ is an isometric embedding of the Hermitian symmetric space
$U/K$ into Euclidean space
$\cu$.    A direct computation shows that the $(a,2)$-flow in the
$U$-hierarchy on $\cs(\R,
\cu_a^\perp)$ is
\refeq[hl]$$u_t= [a,u_{xx}] - {1\over 2} [u, [u, [a,u]]].$$
When $U/K$ is the the complex Grassmannian $\Gr(k, \C^n)$, then $a$
is given by \refhq{}
and equation \refhl{} is the MNLS with $u=\pmatrix{0&q\cr -q^*&0\cr}$.

\ms
It is well-known in soliton theory (cf. [TU2]) that there are two Poisson operators
    on $\cs(\R, \cu_a^\perp)$ so that flows in the $U$-hierarchy on $\cs(\R,
\cu_a^\perp)$ are commuting Hamiltonian flows with respect to both
Poisson structures.  The first
Poisson operator is
\refeq[hb]$$J_a(v)= [v,a].$$
The second Poisson operator $P_u$ is defined as follows:
\refeq[ez]$$P_u(v) = v_x+ \pi_1([u,v]) + [u,h],$$ where $$h(x)=
-\int_{-\infty}^x
\pi_0([u(s),v(s)])ds,$$
and $\pi_0, \pi_1$ are the orthogonal projection of $\cu$ onto $\cu_a$
and $\cu_a^\perp$ respectively.
Moreover, for each $k$,
\refeq[ho]$$(J_k)_u(v)= J_a(J_a^{-1}P_u)^k= (P_uJ_a^{-1})^kJ_a$$ is
also a Poisson
operator on
$\cs(\R,\cu_a^\perp)$.

The Hamiltonian of the $(b,j)$-flow with respect to $J_a$ is  $F_{b,j}:\cs(\R,
\cu_a^\perp)\to \R$ defined by
\refeq[ha]$$F_{b,j}(u)= -{1\over j+1}\int_{-\infty}^\infty \li
Q_{b,j+2}(u), a\ri \ dx.$$
In other words
\refeq[ig]$$\K F_{b,j}(u)=\pi_1( Q_{b,j+1}(u)),$$
and the $(b,j)$-flow is
$$u_t= J_a(\K F_{b,j}(u)).$$
It follows from the definition of $P_u$, the recursive formula
\refad{}, and \refig{}  that
the $(b,j)$-flow can also be written as
\refeq[gi]$$u_t=(J_k)_u(\K F_{b,j-k}(u)), \quad k\leq j.$$
So the $(b,j)$-flow is Hamiltonian with respect to $J_k$
for all $k\leq j$.

\ms
Next we explain how to construct geometric curve flows on Adjoint
orbits from flows in
the $U$-hierarchy.  Let $u$ be a solution of the $(b,j)$-flow in the
$U$-hierarchy on
$\cs(\R,\cu_a^\perp)$. Since its Lax pair $\o_\l$ defined by \refex{}
is flat for all $\l$,
$\o_0= udx + Q_{b,j}(u) dt$ is a flat $\cu$-valued connection 1-form.
Hence there exists
$g:\R^2\to U$ such that
\refeq[ik]$$g^{-1}g_x= u, \quad g^{-1}g_t= Q_{b,j}(u).$$
Set
\refeq[ii]$$\g(x,t)= g(x,t)\ a\  g(x,t)^{-1}.$$
Then $\g(\cdot, t)$ is a family of curves on the Adjoint $U$-orbit
$M_a$ in $\cu$.  Take the
$t$ derivative of \refii{} to get
\refeq[ij]$$\g_t= g[Q_{b,j}(u),a]g^{-1}.$$
We will show that when $M_a$ is a Hermitian symmetric space and
$(b,j)=(a,2)$, the
right hand side of \refij{} is equal to
$$[\g,\g_{xx}] =[\g, \K_{\g_x}\g_x].$$
This implies that solutions of MNLS and \refhl{} give rise to solutions of
the Schr\"odinger flow \refhr{} on
$\Gr(k,\C^n)$ and Hermitian symmetric space respectively.
For $j=1$, the $(b,1)$ flow is the $n$-wave equation, and the right hand side
of  \refij{} is equal to
$$(gbg^{-1})_x= (\hat b(\g))_x= -A_{\hat b(\g)}(\g_x).$$   So solutions of
the $n$-wave equation  give rise to  solutions of the
Ferapontov  flow \refel{}. 
However, some natural  questions come up in this
construction:

\ms

\item {(1)} Since the solution $g$ of \refik{} is only unique up
to left  multiplication, the
corresponding solution $\g$ of \refij{} is unique up to conjugation.
Can we normalize the
curves in
$M_a$ so that the correspondence between $u$ and $\g$ is unique and  preserves
the flow?
\item {(2)} Is \refij{} a geometric curve flow on $M_a$? In other
word, can the right hand side of \refij{} be written as some geometric quantity
$H_{b,j}(\g)$?
\item {(3)} Can the procedure of constructing solutions of \refij{}
from solutions of the
$(b,j)$-flow be reversed?
\item {(4)} Is there a Hamiltonian formulation of the flow \refij{}?
If yes, what is its relation to
the Hamiltonian theory of the $(b,j)$-flow?

\ms

\ni When the $(b,j)$-flow is the MNLS, Terng and Uhlenbeck chose a
normalization for curves
on $M_a$ at $-\infty$ and were able to answer the above questions in
a satisfactory way
([TU2]).  In this paper we generalize their results to any flows in the
$U$-hierarchy.  To
do this, we need to recall the development map constructed in [TU2] next.

Let $M_a$ be the Adjoint $U$-orbit at $a$ in $\cu$,   and $C_a(\R,M_a)$
the space of all smooth curves $\g:\R\to M_a$ such that
$\lim_{x\to -\infty} \g(x)=a$ and  $\g_x\in \cs(\R, \cu)$.
The following results were proved in [TU2]:
\item {(i)} Given $\g\in C_a(\R,M_a)$, there exists a unique
$g:\R\to U$ such that
$$\g= gag^{-1}, \quad\lim_{x\to -\infty}g(x)=e, \quad g_x\in
\cs(\R,\cu_a^\perp),$$
where $e\in U$ is the identity element.
\item {(ii)} The {\it development map\/}
$\Phi:C_a(\R,M_a)\to
\cs(\R,\cu_a^\perp)$ defined by $\Phi(\g)= g^{-1}g_x$
is a bijection.
\item {(iii)} Since $\Phi$ is a bijection,
$\Phi^*(J_k)$ is a Poisson operator on $C_a(\R,M_a)$.  Moreover,
$\Phi^*(J_2)$ is a zero
order Poisson operator, and
$\Phi^*(J_k)$ ($k\geq 3$) is a order
$(k-2)$ non-local Poisson operator. Here a non-local Poisson operator
is said to have order
$k$ if it involves derivatives  up to order
$k$ and antiderivatives.   Note that although
$J_k$ has order
$k$, the pullback $\Phi^*(J_k)$ has order $k-2$.

\ms
Next we give a geometric interpretation of $\Phi^*(J_2)$
and
$\Phi^*(J_3)$. There is a natural zero
order Poisson operators on $C_a(\R,M_a)$ obtained by identifying the
Adjoint $U$-orbit
$M_a$ as the coadjoint orbit at $\ell_a\in
\cu^*$, where $\ell_a(x)= \li x, a\ri$.  So $M_a$ is equipped with
the coadjoint orbit
symplectic form and the corresponding Poisson operator at $y\in M_a$
is $-\ad(y)$. Hence
it induces a natural zero order Poisson operator on $C_a(\R,M_a)$:
\refeq[fl]$$J_\g(v)= [v,\g].$$
In fact,  $J= \Phi^*(J_2)$.

\ms
The Poisson operator $\L= \Phi^*(J_3)$, which is non-local and of
first order,  can be
described geometrically as follows:  Given $v\in T(C_a(\R,M_a))_\g$,
there exists a unique
vector field $\eta$ along $\g$ normal to $M$ such that $(v+\eta)_x$
is tangential.  Then the Poisson operator $\L= \Phi^*(J_3)$ is given by
$\L_\g(v)= (v+\eta)_x$.   Moreover,
$$\Phi^*(J_k)_\g=
(\L_\g J_\g^{-1})^{k-2}J_\g.$$

When $a$ is regular, the Poisson operator $\L$ is the same as the
Poisson
operator constructed by Ferapontov in [F1] using the Dirac reduction of
the Poisson
operator
$d_x$ on $C(\R, \cu)$ to
$C(\R, M_a)$.

\ms
Below are some of our results:
\ms
\item {(i)} the curve flow on $M_a$ corresponding to the
$(b,j)$-flow  under the development
map $\Phi$ is
\refeq[ep]$$\g_t= - (\L_\g J_\g^{-1})^{j-1}(A_{\hat
b(\g)}(\g_x))=(\L_\g J_\g^{-1})^{j-1}(\K H_b(\g)),$$
where  $$H_b(\g)=\int_{-\infty}^\infty \li \g(x),b\ri dx.$$
\item {(ii)} The Poisson operators $\L$, $J$ and the constant of
motions  $F_{b,j}\circ
\Phi$ for  \refep{} can be expressed in geometric terms.
\item {(iii)} The curve flow \refep{} is Hamiltonian with respect to
$J$ and $\L$, and is
completely integrable.
\item {(iv)} If $M_a$ is isometric to a Hermitian symmetric space, then the
curve flow on $M_a$ corresponding to the $(a,2)$-flow \refhl{} under
$\Phi$ is the
Schr\"odinger flow on $M_a$.
\item {(v)} If $M_a$ is a principal Adjoint $U$-orbit in $\cu$ and
$f:\cu\to \R$ is a polynomial invariant under the Adjoint action,
then the curve flow
\refeq[ir]$$\g_t= (\K f(\g))_x$$
with constraint $\g(x,t)\in M_a$, is the Ferapontov flow \refel{} on
$M_a$ with $b=\K
f(a)$.  Moreover, the collection of curve flows
$$\g_t = (\L_\g J_\g^{-1})^{j-1}(\K h(\g)),$$
with $h:\cu\to \R$ a polynomial invariant under the Adjoint
$U$-action and $j\geq 1$, is a
hierarchy of commuting Hamiltonian flows on $C_a(\R,M_a)$ with respect to both
Poisson operators $J$ and $\L$.
\ms

Next we explain how to construct the hierarchy of commuting flows associated
to a symmetric space.    Let
$\s$ be the involution on
$U$ such that
$K$ is the fixed point of $\s$, and
$\cp$ the $-1$ eigenspace of $d\s_e$.  Then $U/K$ is a symmetric
space.  It is known
that (cf.~[TU1]):
\item {(i)} The subspace $\cs(\R,\cu_a^\perp\cap \ck)$ is invariant
under all $(b,j)$-flows
with odd $j$.  The collection of these restricted flows is
called the {\it
$U/K$-hierarchy\/}.
\item {(ii)} If $k$ is odd, then the restriction of the  Poisson
operator
$J_k$ to  the invariant
submanifold
$\cs(\R, \cu_a^\perp\cap \ck)$  is
again a Poisson operator. The odd flows in the
$U/K$-hierarchy are Hamiltonian with respect to these Poisson structures.

\ms
Let $a\in \cp$, and $M_a$ and $N_a$ denote the Adjoint $U$-orbit in
$\cu$ and Adjoint $K$-orbit in $\cp$ at $a$ respectively.  So
$N_a\subset M_a$, and
$C(\R,N_a)$ is a submanifold of $C_a(\R,M_a)$.   We will show that
$C_a(\R,N_a)$ is invariant under the curve flow \refep{} if $j$ is
odd.  Moreover, if $k$ is
odd, then
$\Phi^*(J_k)$ induces a Poisson operator on $C_a(\R, N_a)$.

\ms

This paper is organized as follows:    We write down
     curve flows on Adjoint $U$-orbits corresponding to flows in the
$U$-hierarchy under the
development map as geometric flows, and express the corresponding
Poisson operators and
constant of motions in geometric terms  in section 2.  We construct B\"acklund
transformations and finite type solutions of these geometric curve
flows in section 3, and
consider the curve flows corresponding to flows in the
$U/K$-hierarchy in section 4.   Finally,
we study the curve flow
\refel{} on a principal orbit of the isotropy representation of a
symmetric space as a
hydrodynamic system in section 5.
\bs

\newsection Integrable curve flows on Adjoint orbits.\par

Let $M_a$ denote the Adjoint $U$-orbit at $a$ in $\cu$.
In this section, we write down the curve flows, Poisson structures,
and commuting
Hamiltonians on $C_a(\R,M_a)$ corresponding to soliton
flows in the
$U$-hierarchy via the development map in geometric terms.

    If $f:X\to Y$ is a diffeomorphism and $w$ is a symplectic form on
$Y$, then the pull back
$f^*(w)$ is a symplectic form on $X$ and $f:(X,f^*(w))\to (Y,w)$ is a
symplectic diffeomorphism.
If $g$ and $h$
are Riemannian metrics on
$X$ and
$Y$ respectively, then there exists section $B$ of $L(TX,TX)$ that
relates the metrics
$g$ and
$f^*(g)$ on
$X$ as follows:
$$f^*(g)_x(v_1,v_2)= h_{f(x)}(df_x(v_1), df_x(v_2))= g_x(B_x(v_1), v_2).$$
Let $J$ be the Poisson operator corresponding to $w$ on $Y$, i.e.,
$$w_y(v_1, v_2)= h_y(J_y^{-1}(v_1), v_2).$$
    A direct computation shows that the Poisson operator $f^*(J)$
corresponding to
$f^*(w)$ is
\refeq[gb]$$f^*(J)_x = df_x^{-1}\circ J_{f(x)} \circ df_x \circ B_x^{-1}.$$
The gradient $\K H$ of $H:Y\to \R$ is defined by
$$dH_y(v)= h_y(\K H(y), v).$$
The Hamiltonian flow for $H$ on $Y$ with respect to $w$ is
\refeq[gc]$${dy\over dt}= J_{y(t)}(\K H(y(t))).$$
The Hamiltonian equation for $H\circ f:X\to \R$ with respect to $f^*(J)$  is
\refeq[gd]$${dx\over dt}= f^*(J)_x(\K (H\circ f)(x)) =
df_x^{-1}(J_{f(x)}(\K H (f(x))).$$
It is clear that the
Hamiltonian flow for $H\circ f$ on $X$ maps to the Hamiltonian flow
for $H$ on $Y$.
Equations \refgc{} and \refgd{} can be viewed as the same equation
written in different
coordinate systems.

\ms

Let $\Phi:C_a(\R,M_a)\to \cs(\R,\cu_a^\perp)$ be the development map
on the Adjoint
$U$-orbit $M_a$ given in section 1.  To compute $\Phi^*(J_k)$, we
need to specify the
metrics.
    The space $\cs(\R, \cu_a^\perp)$ is
equipped with the $L^2$ inner product
$$(u_1, u_2)= \int_{-\infty}^\infty \li u_1(x), u_2(x)\ri \ dx.$$
The tangent space of $\cc_a(\R, M_a)$ at $\g$ is the space of vector
fields $\xi$ tangent to
$M_a$ along $\g$ and $\xi_x\in \cs(\R,\cu)$.   Let
$ds^2$  denote the $L^2$
metric on $\cc_a(\R, M_a)$ defined by
$$ds^2_\g( \xi_1, \xi_2)= \int_{-\infty}^\infty \li \xi_1(x),
\xi_2(x)\ri \ dx$$ for
tangent fields  $\xi_1, \xi_2$ along $\g$.

The following is proved in [TU2]:

\refclaim[ge] Proposition.  Let $\xi\in T(C_a(\R,M_a))_\g$. Then:
\item {(i)} There exist $g:\R\to U$ and
$v:\R\to \cu_a^\perp$ so that $\lim_{x\to -\infty} g(x)=e$,
$\g=gag^{-1}$, $u=g^{-1}g_x= \Phi(\g)$, and $\xi=gvg^{-1}$.
\item {(ii)} There exists a unique vector field $\eta$ along $\g$
such that $\eta(x)$ is
normal to $M_a$, $\lim_{x\to -\infty}\eta(x) = 0$, and
$(\xi+\eta)_x$ is tangent
to  $M_a$.
\item {(iii)}
\refeq[hm]$$d\Phi_\g(\xi) = P_uJ_a^{-1}(v),$$
where $J_a$ and $P_u$ are the Poisson operators on $\cs(\R, \cu_a^\perp)$
defined by \refhb{} and \refez{} respectively.
\item {(iv)} Let $J_3=(P_uJ_a^{-1})^3J_a$ be the Poisson operator
defined by \refho{}, and
$\L=\Phi^*(J_3)$. Then
\refeq[hf]$$\L_\g(\xi)
=g(P_u(v))g^{-1}= (\xi+\eta)_x,$$  where $\eta$ and $g$ are given in (ii).\ei

A direct computation implies

\refclaim[hn] Proposition ([TU2]).  Let $J_2$ be the Poisson operator on
$\cs(\R,\cu_a^\perp)$ defined by $(J_2)_u= P_uJ_a^{-1}P_u$, $J$ the natural
Poisson operator on $C_a(\R,M_a)$ given by \reffl{}, and $\Phi$ the development
map.    Then
$\Phi^*(J_2)= J$.

When $M_a$ is a principal Adjoint $U$-orbit, we show below that the
Poisson operator $\L$
can be written in terms of geometric invariants of $M_a$ as a
submanifold of the Euclidean
space $\cu$.

\refclaim[gk] Theorem.  Let $\ct$ be a maximal abelian
subalgebra of $\cu$, $\{a_1, \ldots, a_k\}$ an orthonormal basis of
$\ct$, and
$a=a_1$  regular.
Let $M_a$ be the Adjoint $U$-orbit at $a$ in $\cu$, and $\L=\Phi^*(J_3)$  the
induced Poisson operator on $C_a(\R,M_a)$.  Then
\refeq[fm]$$\L_\g(\xi) =\K_{\g_x}\xi -\sum_{i=1}^{k} h_i A_{\hat
a_i(\g)}(\g_x),$$
where $h_i(x) = -\int_{-\infty}^x \li \II(\xi(s), \g_x(s)),
\hat a_i(\g(s))\ri ds$ for $1\leq \a\leq  k$, $\hat a_i$ is the
parallel normal field  on
$M_a$ defined by \refeo{}, and  $\II$ is the second fundamental form
of $M_a$ as a
submanifold of
$\cu$.

\proof
Let $u=\Phi(\g)$, and $\xi, v, g$ as in Proposition \refge{}.  So we have
$$\g=gag^{-1}, \ \ u=g^{-1}g_x, \ \, \xi= gvg^{-1}.$$
     By Proposition \refge{} (iv) and the
formula \refez{} for $P_u$ , we have
$$\L_\g(\xi)= gP_u(v)g^{-1}= g(v_x+[u,v]_1+[u,h])g^{-1},$$ where
$h_x=-[u,v]_0$, and
$[u,v]_0$ and
$[u,v]_1$ denote the projection of $[u,v]$
onto $\cu_a=\ct$ and $\cu_a^\perp=\ct^\perp$ respectively.

Given $\zeta\in \cu$, let $p_x(\zeta)$ and $p_x^\perp(\zeta)$ denote
the orthogonal
projection of $\zeta$ onto $T(M_a)_x$ and $\nu(M_a)_x$ respectively.
Let $\K$ denote
the Levi-Civita connection of the induced metric on $M_a$.  Then
$$\eqalign{(\K_{\g_x}\xi )(x) &= p_{\g(x)}((gvg^{-1})_x) =
p_{\g(x)}(g(v_x+[u,v])g^{-1})\cr & = g(v_x+ [u,v])_1 g^{-1} = g(v_x +
[u,v]_1)g^{-1}.\cr}$$

    But
$$(ghg^{-1})_x= g(h_x+ [u,h])g^{-1}.$$ Since $h(x)\in \ct$ and $a_1,
\ldots, a_k$
is a basis of $\ct$, there exist $h_1, \ldots, h_k$ such that
$$h(x)= \sum_{i=1}^k h_i (x)a_i.$$ Therefore
$$gh_xg^{-1}= g\left(\sum_{i=1}^k (h_i)_x a_i\right)g^{-1} =
\sum_{i=1}^k (h_i)_x
\hat a_i(\g(x)).$$
Let $\K^\perp$ denote the induced normal connection on $\nu(M)$.
Then
$$\K^\perp_{\g_x} (ghg^{-1})= gh_xg^{-1}=\sum_{i=1}^n (h_i)_x \hat a_i.$$
Definition of the second fundamental form implies that
$$\eqalign{\II(\xi, \g_x) &= p_{\g(x)}^\perp((gvg^{-1})_x) = g[u,v]_0g^{-1}= -
gh_xg^{-1} \cr & = - \sum_{i=1}^k (h_i)_x \hat a_i \cr &=
\sum_{i=1}^k\li \II(\xi,
\g_x),\hat a_i\ri \hat a_i.\cr}$$
So $(h_i)_x = -\li\II(\xi,\g_x), \hat a_i(\g(x))\ri$. But $\lim_{x\to
-\infty}\eta(x)=0$. So $h_i$ is given by the formula in the Proposition. \qed

\refpar[] Remark. Let $M^n$ be a submanifold of $\R^{n+k}$ with flat
normal bundle, and
$(e_{n+1}, \ldots, e_{n+k})$ a parallel orthonormal normal frame field on
$M$.  Note that $\xi$ is a tangent vector of $C(\R,M)$ at $\g$ iff
$\xi$ is a vector field
along $\g$ tangent to $M$.   Let $Z_\g$ be the operator on
$TC(\R,M)_\g$ defined by
\refeq[fm]$$\eqalign{Z_\g(\xi)& =\K_{\g_x}\xi +\sum_{\a=n+1}^{n+k} h_\a
A_{\hat e_\a(\g)}(\g_x),
\quad {\rm where\/}\cr  h_\a(x) &= \int_{-\infty}^x \li \II(\xi(s), \g_x(s)),
e_\a(\g(s))\ri ds
\quad {\rm for\ } n+1\leq \a\leq n+ k,\cr}$$
and $\K$ is the Levi-Civita connection of the induced metric on $M$.
Ferapontov proved in
[F1] that $Z$ is a Poisson operator on $C(\R,M)$.  The normal bundle of
a principal Adjoint
$U$-orbit
$M_a$ is flat.  Proposition
\refgk{} proves that
$\Phi^*(J_3)=\L = Z$.

    \refclaim[gg] Theorem ([TU2]).  Let $J_k$ be the  Poisson structures
on $\cs(\R,\cu_a^\perp)$ defined by  \refho{}, and
$L_{k-2}= \Phi^*(J_k)$  the induced Poisson
structures on $C_a(\R,M_a)$.  Then
\refeq[eu]$$(L_j)_\g=J_\g ( J_\g^{-1}\L_\g)^j= (\L_\g J_\g^{-1})^jJ_\g.$$

It follows from \refgd{} that the flow on $C_a(\R,M_a)$
corresponding to the
$(b,j)$-flow in the $U$-hierarchy on $\cs(\R, \cu_a^\perp)$ is
\refeq[gh]$$\g_t= d\Phi_\g^{-1}([Q_{b,j+1}(u),a]).$$
Propositions \refge{} and \refgi{}
imply that

\refclaim[hs] Proposition. The curve flow on $M_a$ corresponding to the
$(b,j)$-flow in the
$U$-hierarchy on $\cs(\R, \cu_a^\perp)$ is
\refeq[gl]$$\g_t= g[Q_{b,j}(u),a]g^{-1}=[gQ_{b,j}(u)g^{-1}, \g],$$
where $\g=gag^{-1}$, $\lim_{x\to -\infty} g(x,t)= e$, and
$u=g^{-1}g_x= \Phi(\g)$.
    Moreover, \refgl{}  is the Hamiltonian flow for
$F_{b,j-k-2}\circ \Phi$ with respect to $L_k=\Phi^*(J_{k+2})$.

Next we want to express \refgl{} and $F_{b,j}\circ \Phi$
in geometric terms.  This follows from

\refclaim[ht] Proposition.  Let $\g, u, g$ be as in Proposition \refhs{}.  Then
$gQ_{b,j}(u)g^{-1}$ can be expressed in
terms of geometric invariants of
$M_a$ along $\g$.

\proof
We prove this proposition by induction.
The recursive formula
\refad{} implies $$Q_{b,1}(u)=\ad(b)\ad(a)^{-1}(u).$$  A direct
computation gives
$$\eqalign{& \g_x= g[u,a]g^{-1}, \cr
& gug^{-1}= J_\g^{-1}(\g_x),\cr
& gQ_{b,1}(u)g^{-1}= J_\g^{-1}(g[Q_{b,1}(u),a]g^{-1}) =
J_\g^{-1}(g[u,b]g^{-1}).\cr}$$
The parallel normal field $\hat b$ defined by \refeo{} is $\hat
b(gag^{-1})= gbg^{-1}$.
The shape operator for
$M_a$ along $\hat b$ is
$$A_{\hat (\g)}(\g_x) = -(\hat b(\g))_x = -(gbg^{-1})_x = -g[u,b]g^{-1}.$$
Hence
\refeq[it]$$gQ_{b,1}(u)g^{-1}= -J_\g^{-1}(A_{\hat b(\g)}(\g_x))=
J_\g^{-1}(\hat
b(\g)_x).$$
Suppose we have expressed
$gQ_{b,j}(u)g^{-1}$ geometrically. Let
$\pi_0, \pi_1$ be the
orthogonal projection onto $\cu_a$ and $\cu_a^\perp$ respectively,
and $p_x$ and
$p_x^\perp$ the orthogonal projection of $\cu$ onto $\cu_x=\nu(M_a)_x$ and
$\cu_x^\perp=T(M_a)_x$ respectively.   The recursive formula
\refad{} and Proposition \refge{}  imply that
\refeq[iu]$$g\pi_1(Q_{b,j+1}(u))g^{-1}=
J_\g^{-1}\L_\g(g\pi_1(Q_{b,j}(u))g^{-1}),$$ and
\refeq[hx]$$\eqalign{(g\pi_0(Q_{b,j+1}(u))g^{-1})_x&=
-g\pi_0([u,\pi_1(Q_{b,j+1}(u))])g^{-1}\cr &=
-p_{gag^{-1}}^\perp([gug^{-1},
g\pi_1(Q_{b,j+1}(u))g^{-1}])\cr &= -p_{gag^{-1}}^\perp([J_\g^{-1}(\g_x),
g\pi_1(Q_{b,j+1}(u))g^{-1}]).\cr}$$
\ni The induction hypothesis and \refiu{} imply that
$g\pi_1(Q_{b,j+1}(u))g^{-1}$ can be expressed geometrically. Then
\refhx{} implies
that $g\pi_0(Q_{b,j+1}(u))g^{-1}$ can also expressed as geometric terms.
\qed

    Use \refit{}, \refiu{} to get
\refeq[hy]$$g[Q_{b,j}(u),a]g^{-1}= (\L_\g J_\g^{-1})^{j-1}((\hat
b(\g))_x) = -(\L_\g
J_\g^{-1})^{j-1}(A_{\hat b(\g)}(\g_x)).$$
   The next theorem follows from \refhy{} and Proposition \refhs{}.

\refclaim[gm] Theorem.  Suppose $a\in \cu$ is regular, and $b\in
\cu_a$. Then the curve
flow  on $C_a(\R,M_a)$ corresponding to the $(b,1)$-flow in the
$U$-hierarchy  is
\refeq[gn]$$\g_t=(\hat b\circ \g)_x= -A_{\hat b(\g)}(\g_x),$$
where $A_{\hat b(\g)}$ is the shape operator along $\hat b(\g)$.
In general,  the curve flow on $C_a(\R,M_a)$ corresponding to the
$(b,j)$-flow is
\refeq[ep]$$\g_t= (\L_\g J_\g^{-1})^{j-1}(b(\g))_x=
-(\L_\g J_\g^{-1})^{j-1}(A_{\hat b(\g)} (\g_x)).$$
Moreover,
\item {(i)} $F_{b,j}\circ \Phi$ can be expressed in terms of shape
operators of $M_a$
and the Poisson operators $J$ and $\L$,
\item {(ii)}  the flow \refep{} is the Hamiltonian flow of
$F_{b,j-k}\circ
\Phi$ on
$C_a(\R,M_a)$ with respect to $L_{k-2}=\Phi^*(J_k)$, where $F_{b,j}$ is the
functional defined by \refha{} and $\Phi$ is the development map.
\ei

\refpar[] Example. Let $a\in \cu$ be a regular element.  The curve
flow \refep{} on the
Adjoint orbit $M_a$ corresponding to the $(b,2)$-flow is
$$\g_t= -\big(J_\g^{-1}(A_{\hat b(\g)}(\g_x))\big)_x.$$
\ms

When $b=a$ (without the assumption that $a$ is regular), we have

\refclaim[go] Corollary.  For $a\in \cu$,  the curve flow on
$C_a(\R,M_a)$ corresponding to the $(a,j)$-flow in the $U$-hierarchy is
\refeq[hu]$$\g_t= (\L_\g J_\g^{-1})^{j-1}(\g_x).$$

Since \refep{} is the flow corresponding to the
$(b,j)$-flow under the  development map
$\Phi$ and $u$ is a solution of the $(b,j)$-flow,
we get that
$\g(\cdot,t)=\Phi^{-1}(u(\cdot, t))$ is a solution of \refep{}.
Proposition
\refhs{} states that
\refeq[gl]$$\g_t= g[Q_{b,j}(u), a]g^{-1},$$
where $g$ is the solution of
\refeq[gv]$$g^{-1}g_x= u, \quad \lim_{x\to -\infty} g(x,t)=e.$$
On the other hand, we have seen in
section 1 that if $k:\R^2\to U$ solves
\refeq[ik]$$k^{-1} k_x= u, \quad k^{-1}k_t= Q_{b,j}(u),$$
then
$$\ti \g = k ak^{-1}$$ also satisfies \refgl{}.  Next we show that
$g$ also satisfies
\refik{}.

\refclaim[gs] Proposition.  Let $u$ be a solution of the $(b,j)$-flow
on $\cs(\R,
\cu_a^\perp)$, $\o_\l$ the corresponding Lax pair \refex{},  and $k$
the trivialization of
the Lax pair $\o_0$, i.e.,
$$k^{-1}dk= u\ dx + Q_{b,j}(u)\ dt, \quad k(0,0)=e.$$
Then there exists a constant $c\in U$ such that $\lim_{x\to -\infty}
k(x,t)= c$,
$$\g= c^{-1}kak^{-1}c$$
is a solution of  the curve flow \refep{} on $C_a(\R,M_a)$, and
$\Phi(\g(\cdot, t))= u(\cdot,
t)$.

\proof
For fix $t$, let $g(\cdot, t)$ be the map from $\R$ to $U$ satisfying \refgv{}.
Claim that $g^{-1}g_t= Q_{b,j}(u)$.   Too see this, we note that because both
$g$ and $k$ satisfy
$g^{-1}g_x= k^{-1}k_x= u$, there exists
$C(t)$ such that
$$g(x,t)=C(t) k(x,t).$$ Compute directly to get
\refeq[fa]$$g^{-1}g_t= k^{-1}k_t + k^{-1}C^{-1}C_t k = Q_{b,j}(u) +
k^{-1}C^{-1}C_t
k.$$
Since $\lim_{x\to -\infty} g(x,t)= e$, $\lim_{x\to -\infty}
g_t(x,t)= 0$.
Recall that $Q_{b,j}(u)\in
\cs(\R, \cu)$  if $u(\cdot, t)\in \cs(\R,\cu_a^\perp)$ (\refjd{}).  Hence
$$\lim_{x\to -\infty} Q_{b,j}(u)(x,t)=0.$$
By \reffa{}, we conclude
\refeq[fw]$$\lim_{x\to -\infty} k^{-1}C^{-1}C_t k=0.$$
But
\refeq[fx]$$\N k^{-1}C^{-1}C_tk\N =\N C^{-1}C_t\N.$$
     It follows from \reffw{} and \reffx{} that $C^{-1}C_t=0$. Hence
$C(t)$  is a
constant, $g=Ck$, and we have proved the claim.  \qed

\refclaim[gw]  Corollary.  Let $u$ be a solution of the $(b,j)$
flow, $k, c$ as in
Proposition \refgs{}, and $g=c^{-1}k$.   Then the gauge
transformation of the Lax pair
$\o_\l$ (\refex{}) by $g$,
$$\g \l \ dx + (gbg^{-1}\l^j + gQ_{b,1}(u)g^{-1} \l^{j-1}+ \cdots +
gQ_{b,j-1}(u)g^{-1}\l)\ dt$$
    is a Lax pair of the curve flow \refep{}.

\refpar[il] Example.  The curve flow \refgn{} on $M_a$
  has a Lax pair
$$\g\l \ dx + \hat b(\g) \l \ dt.$$

\ss
\refpar[im] Example ([F4]).  Let $M_a$ be a principal Adjoint
$U(n)$-orbit in $u(n)$.  If $k$ is
a natural number, then
$b=i^{k-1}a^k\in u(n)_a$, and the curve flow
\refel{} on
$M_a$ corresponding to the
$(b,1)$-flow becomes
$$\g_t= i^{k-1}(\g^k)_x$$
with $\g(x,t)\in M_a$.
Its Lax pair is $\g\l \ dx + i^{k-1} \g^k \l \ dt$.

\ss
\refpar[in] Example.  Let $a\in \cu$ be regular, and $b\in \cu_a$.
It is known (cf.~[Te1]) that
there exists a polynomial $f:\cu\to \R$ invariant under the Adjoint
action such that $\K f(a)=
b$. Since $f$ is Ad-invariant, $\K f$ is equivariant.  Hence $\K f\n
M_a$ is a parallel normal
field of $M_a$.  So the Ferapontov flow  \refel{} for $b= \K f(a)$ becomes
\refeq[hv]$$\g_t= (\K f(\g))_x, \quad \g(x,t)\in M_a.$$
Its Lax pair is $\g\l \ dx + \K f(\g) \l \ dt$.
Since the $(b,1)$-flow \refab{} is a completely integrable
Hamiltonian system with respect to
the Poisson operator $J_2$ and $J_3$, and
\refhv{} is the flow corresponding to the
$(b,1)$-flow under the development map $\Phi$, the curve flow
\refhv{} is completely
integrable with respect to the Poisson operators $J$ and
$\L=\Phi^*(J_3)$.  Moreover, the
flows
$$\g_t= (\L_\g J_\g^{-1})^j((\K h(\g))_x),$$
with $h:\cu\to \R$ a polynomial invariant under the Adjoint
$U$-action and $j\geq 0$ an
integer, are commuting Hamiltonian flows on $C_a(\R,M_a)$.

\ss

\refpar[io] Example.  Let $U/K$ be a compact Hermitian symmetric
space. Then there
exists $a\in \ck$ such that $\cu_a=\ck$ and $\ad(a)^2=-\Id$ on
$\cp=\cu_a^\perp$.
Moreover, the Adjoint $U$-orbit $M_a$ at $a$ in $\cu$ is an isometric
embedding of
the Hermitian symmetric space $U/K$. In fact,  the induced metric on
$M_a$ is the standard K\"ahler metric on $U/K$,  $j_\g(\xi)=[\g,\xi]$
is  the complex
structure on $U/K$.  Use condition \refad{} and \refjd{} and a direct
computation to see that the $(a,2)$-flow in
the
$U$-hierarchy on
$\cs(\R, \cu_a^\perp)$ is \refhl{}:
$$u_t= [a,u_{xx}] - {1\over 2} [u, [u, [a,u]]].$$
The corresponding curve flow
\refhu{} on $M_a$ is
\refeq[fh]$$\g_t= \L_\g J_\g^{-1}(\g_x).$$
Claim that this flow is the
Schr\"odinger flow on
the Hermitian symmetric space $U/K$:
\refeq[ib]$$\g_t= j_\g(\K_{\g_x}\g_x) =[\g, \K_{\g_x}\g_x].$$
     To see this, let $g$ be the solution of
$$g^{-1}g_x=u, \quad \lim_{x\to -\infty}
g(x,t)=e.$$  Set $\g= gag^{-1}$.  Then $$\g_x= g[u,a]g^{-1}.$$
We compute the right hand side of
\reffh{} as follows:
$$\L_\g J_\g^{-1}(\g_x)= \L_\g(J_\g^{-1}(g[u,a]g^{-1}))
=\L_\g(gug^{-1}).$$
  It follows from \refez{} and \refhf{} that
$$\L_\g(gug^{-1}) = gu_xg^{-1}.$$
  But
\refeq[ic]$$\K_{\g_x}\g_x= p(g[u,a]g^{-1})_x= p(g[u_x,a]g^{-1} +
g[u,[u,a]]g^{-1}),$$
where $p(v)$ denotes the orthogonal projection of $v$ onto $TM_a$ and
$\K$ is the Levi-Civita connection of the induced metric on
$M_a$.  But
$$(TM_a)_\g= g\cu_a^\perp g^{-1}=g\cp g^{-1}, \quad \nu(M_a)_\g= g\cu_a
g^{-1}= g\ck g^{-1}.$$ Since $M_a= U/U_a = U/K$ is a symmetric space,
$$[\ck,\ck]\subset \ck,
\quad [\ck,\cp]\subset \cp, \quad [\cp, \cp]\subset \ck.$$
So $[\cp, [\cp,\ck]]\subset \ck$, which implies that
$g[u,[u,a]]g^{-1}$ is normal to
$M_a$.  Hence
$$\K_{\g_x}\g_x = g[u_x,a]g^{-1}.$$
   So we have
$$j_\g(\K_{\g_x} \g_x)= [\g, g[u_x,a]g^{-1}] = g[a, [u_x,a]]g^{-1} =
gu_xg^{-1}.$$
Here we use the fact that $-\ad(a)^2=-\id$.
Therefore we have proved the claim, i.e., \reffh{} is the
Schr\"odinger flow \refib{} on
$M_a$.

Next we claim that the right hand side of \refib{} is also equal to
$[\g, \g_{xx}]$. This is done by a  direct computation:
Since $\g_x= g[u,a]g^{-1}$,
$$\g_{xx}=(g[u,a]g^{-1})_x= g([u,[u,a]] + [u_x,a])g^{-1}.$$
Hence
$$[\g,\g_{xx}]= g([a,[u,[u,a]]] + [a,[u_x,a]])g^{-1}.$$
Since $\ad(a)^2= -\Id$ on $\cu_a^\perp$, $[a,[u_x,a]]= u_x$.
Jacobi-identity implies that
$$[a,[u,[u,a]]]= -[u,[[u,a],a]]-[[u,a], [a,u]]= -[u,-u]+0 =0.$$
So $[\g,\g_{xx}]= gu_xg^{-1}$, and the flow \reffh{} becomes
\refeq[gr]$$\g_t=[\g,\K_{\g_x}\g_x]=  [\g,\g_{xx}], \quad \g(x,t)\in M_a.$$
Equation \refib{} has a Lax pair
$$\g\l \ dx + (\g \l^2 + [\g,\g_x]\l \ dt).$$

\refpar[ioi] Example.
  Let $U/K$ be a compact Hermitian symmetric space and let
$a$ be as in Example \refio{}. Then we know that the $(a,2)$-flow in the
$U$-hierarchy on
$\cs(\R,\cu_a^\perp)$ is \refhl{}
$$u_t= [a,u_{xx}] - {1\over 2} [u, [u, [a,u]]].$$
We would
like to write this equation out explicitly for each irreducible 
compact symmetric
space involving classical groups.
\item{(i)} As was pointed out in the introduction, the equation in case of the
Grassmannian $\Gr(k,\C^n)$ is the MNLS.
\item{(ii)} We consider the Grassmannian $\Gr(2,\R^{n+2})=SO(n+2)/SO(2)\times
SO(n)$. We have $U=SO(n+2)$.
The element
$a$ in $so(n+2)$ has the matrix
$$\pmatrix{ 0&-1\cr 1 & 0\cr}$$
in the upper left corner and elsewhere zeros. Its centralizer $\cu_a$ is
$so(2)\times so(n)$ and $\cu^\perp_a$ is the set of matrices of the form
$$\pmatrix{ 0&0&x_1&\cdots & x_n\cr  0&0&y_1&\cdots & y_n\cr
-x_1&-y_1&0&\cdots&0\cr \cdot&\cdot&\cdot&\cdots&\cdot\cr
\cdot&\cdot&\cdot&\cdots&\cdot\cr\cdot&\cdot&\cdot&\cdots&\cdot\cr
-x_n&-y_n&0&\cdots&0\cr}
$$
where $X=(x_1,\dots,x_n)$, $Y=(y_1,\dots,y_n)\in \R^n$. One can clearly
identify
$\cu^\perp_a$ with  $\C^n$ by mapping the above matrix onto $Z=X+iY$.
Under this identification the complex structure $\ad(a)$ on  $\cu^\perp_a$
corresponds to the usual complex structure on $\C^n$.
The equation \refhl{} becomes the following system of equations
$$
\cases{X_t= -Y_{xx} +(X\cdot Y) X - {1\over 2}(3X\cdot X + Y\cdot Y) Y,\cr
Y_t= X_{xx} + {1\over 2} (X\cdot X + 3 Y\cdot Y) X - (X\cdot Y) Y,\cr}
$$ where $X\cdot Y$ is the standard inner product on $\R^n$.
\item{(iii)}  We now consider the Hermitian symmetric space
$SO(2n)/U(n)$. We have $U=SO(2n)$. The element
$a$ is
$$
a={1\over 2}\pmatrix { 0 & -I_n \cr I_n & 0\cr }
$$
whose centralizer $\cu_a$ is $u(n)$. The
  embedding of $U(n)$ into $SO(2n)$ can be described on
the level of Lie algebras as follows:
$Z=A+iB$ in
$u(n)$ where $A$ and $B$ are real matrices ($A$ skew and $B$ symmetric)
is mapped to
$$
\pmatrix { A & B \cr -B & A \cr}
$$
in $so(2n)$. Hence $\cu_a^\perp$ is the set of matrices
$$
  \pmatrix{ X & Y \cr Y & -X \cr}
$$
where $X$ and $Y$ are both skew.
The space $\cu_a^\perp$ can be identified with $\Lambda^2(\C^n)$ by associating
the matrix
$$
\pmatrix{ X & Y \cr Y & -X \cr}
$$
to the form in $\Lambda^2(\C^n)$ whose matrix with respect to the
canonical basis of $\Lambda^2(\C^n)$ is $X+iY$.   Notice that the complex
structure $\ad(a)$ on $\cu_a^\perp$ and the standard one on $\Lambda^2(\C^n)$
correspond under this identification.
The equation \refhl{} now becomes
$$
\cases{X_t = (-Y)_{xx}+[X,[X,Y]]+2Y^3+YX^2+X^2Y,\cr
          Y_t = X_{xx}+[Y,[X,Y]]-2X^3-XY^2-Y^2X.\cr}
$$
\item{(iv)}  We finally consider the Hermitian symmetric space 
$Sp(n)/U(n)$. The
group
$U$ is therefore $Sp(n)$. Recall that the Lie algebra $sp(n)$ consists of
matrices of the form
$$
\pmatrix { A & -\bar B \cr B & \bar A \cr}
$$
where $A$ and $B$ are complex $n\times n$ matrices, $\bar A^t=-A$ and
$B^t=B$.
The element
$a$ is
$$
a={1\over 2}\pmatrix { 0 & -I_n \cr I_n & 0\cr }.
$$
The stabilizer $\cu_a$ is $u(n)$. The embedding of
$U(n)$ into
$Sp(n)$ can be described on the level of Lie algebras as follows:
$Z=A+iB$ in
$u(n)$, where $A$ and $B$ are real matrices ($A$ skew and $B$ symmetric),
is mapped to
$$
\pmatrix { A & -B \cr B & A \cr}
$$
in $sp(n)$. Notice that this means that $\cu_a=u(n)$ consists of the real
matrices in $sp(n)$. It is therefore clear that $\cu^\perp_a$ consists of the
purely imaginary matrices in $sp(n)$. In other words, $\cu_a^\perp$ is the set
of  matrices of the form
$$
i \pmatrix { X & Y \cr Y & -X \cr}
$$
where $X$ and $Y$ are real symmetric matrices.
The space $\cu^\perp_a$ can be identified with the space $S^2(\C^n)$ 
of symmetric
two-forms on $\C^n$ by associating the matrix
$$
\pmatrix{ X & Y \cr Y & -X \cr}
$$
with the symmetric form  whose matrix with respect to the
canonical basis  is $X+iY$. Notice that the complex
structure $\ad(a)$ on
$\cu^\perp_a$ and
the standard one on $S^2(\C^n)$ correspond under this identification.
The equation \refhl{} becomes the system
$$
\cases{X_t = (-Y)_{xx}-[X,[X,Y]]-2Y^3-YX^2-X^2Y,\cr
          Y_t = X_{xx}-[Y,[X,Y]]+2X^3+XY^2+Y^2X,\cr}
$$
where $X, Y$ are real symmetric $n\times n$ matrices.

\ms
   We recall the following result
proved in [Te2]  concerning the operator $P_u$ (Remark 4.1 of [Te2]):

\refclaim[iv] Proposition.  Let $v\in \cs(\R,\cu_a^\perp)$.  Suppose
there exists a
smooth map  $\ti v:\R\to \cu$ such that $\pi_1(\ti v)=v$, $\lim_{x\to
-\infty} \ti
v(x)=0$, and $(\ti v)_x + [u,\ti v]\in
\cu_a^\perp$.  Then $P_u(v)= (\ti v)_x + [u, \ti v]$.

     Next we prove that the Hamiltonian for the curve
flow \refel{} with respect to $\L$ is the height function.

\refclaim[he] Proposition. Let $M_a$ be the principal Adjoint
$U$-orbit at $a$ in $\cu$,
$b\in \nu(M_a)_a$, and
$\hat b$ the parallel normal field on
$M_a$ defined by \refco{}.  Then the Hamiltonian flow  for
$$H_b(\g)=\int_{-\infty}^\infty \li \g(x), b\ri dx$$ with respect
to  $\L=\Phi^*(J_3)$ is the Ferapontov flow \refgn{}:
$$\g_t= (\hat b(\g))_x.$$

\proof
Given $x\in M_a$, let $p_x$ and $p_x^\perp$ denote the orthogonal
projection of $\cp$
onto $T(M_a)_x$ and $\nu(M_a)_x$ respectively. Recall that
$$\nu(M_a)_{gag^{-1}}=g\cu_a g^{-1}, \quad T(M_a)_{gag^{-1}} = g\cu_a^\perp
g^{-1}.$$  The gradient of
$H_b$ at
$\g$ is the vector field along
$\g$ defined as follows:  Suppose
$$\g=gag^{-1}, \quad\lim_{x\to -\infty} g(x)= e, \quad g^{-1}g_x=u\in
\cu_a^\perp$$ as
before.  Then
$$\K H_b(\g)= p_\g(b)= p_{gag^{-1}}(g(g^{-1}bg)g^{-1})=
g\pi_1(g^{-1}bg)g^{-1},$$ where $\pi_1$ is the projection onto $\cu_a^\perp$.
We use \refhf{} and Proposition \refiv{} to compute $\L(\K H_b(\g))$ next.  Let
$$\ti v= g^{-1}bg - b.$$
Then $\pi_1(\ti v)= \pi_1(g^{-1}bg)$ and $\lim_{x\to -\infty} \ti v(x)=0$.  But
$$(\ti v)_x + [u,\ti v]= (g^{-1}bg-b)_x + [u, g^{-1}bg-b]=-[u,b]
\,\in \cu_a^\perp.$$
By Proposition \refiv{}, we have
$$P_u(\pi_1(g^{-1}bg))= -[u,b].$$
Therefore the
Hamiltonian vector
field for
$H_b$ with respect to
$\L$ is
$$\L_\g (\K H_b(\g)) = gP_u(\pi_1(g^{-1}bg))g^{-1} = g[u,b]g^{-1}= (\hat
b(\g))_x.$$
This completes the proof. \qed

\bs

\newsection  B\"acklund transformations and finite type solutions.\par

It is well-known that the $(b,j)$-flow on $\cs(\R, \cu_a^\perp)$ has B\"acklund
transformations, explicit soliton solutions, and finite type
solutions.  So we can
use the development map $\Phi$ and Proposition
\refgs{} to get similar results for the curve flow \refep{} on the
Adjoint $U$-orbit $M_a$.
We explain the construction for the group $SU(n)$.  Other compact
groups can be worked out
in a similar manner.

Let $\g$ be a solution of the curve flow \refep{} on $M_a$.  Given
$z\in \C$ and a
complex linear subspace $V$ of $\C^n$, we will construct a new
solution $\ti\g$ of
\refep{} associated to $\g, z, V$.   This is done using the
B\"acklund transformation of the
solution $u$ of the $(b,j)$-flow corresponding to
$\g$ under the development map $\Phi$. We give a quick review of the
algorithm next.   Let
$u$ be a solution of the $(b,j)$-flow in the $U$-hierarchy, and
$E(x,t,\l)$ the trivialization of the Lax pair $\o_\l$ given by
\refex{} normalized at the origin,
i.e., $E$ is the solution of
$$\cases{E^{-1}E_x= a\l + u, &\cr
E^{-1}E_t= b\l^j + Q_{b,1}(u)\l^{j-1} + \cdots + Q_{b,j}(u),&\cr
E(0,0,\l)=e.&\cr}$$
We call $E$ the {\it frame \/} of the solution $u$.
Let  $\pi$ be the Hermitian projection of $\C^n$
onto $V$, and $\pi^\perp= I - \pi$.  Let
$$f_{z, \pi}(\l)= I +{\bar z-z \over\l-\bar z}\ \pi^\perp.$$
Then $f$ satisfies the $SU(n)$-reality condition, i.e.,
$$f_{z,\pi}(\bar\l)^*f_{z,\pi}(\l) =I.$$
(Here $X^*= \bar X^t$).
Set
$$\ti V(x,t)= E(x,t,z)^*(V),$$
and $\ti \pi(x,t)$ the Hermitian projection of $\C^n$ onto $\ti V(x,t)$.
Then (cf.~[TU1])
\item {(i)} $\ti u= u+ (z-\bar z) [\ti \pi, a]$ is a solution of the
$(b,j)$-flow,
\item {(ii)} $\ti E(x,t,\l)=  f_{z,\pi}(\l)E(x,t,\l)
f_{z,\ti\pi(x,t)}(\l)^{-1}$ is the frame for
the solution $\ti u$.
\ss
\ni  The transformation
$u\mapsto \ti u$ is a {\it B\"acklund transformation\/}.

\ms

By Proposition \refgs{}, there exists a constant $c\in SU(n)$ such that
$$\lim_{x\to -\infty} \ti E(x,t,0)=c$$
for all $t$.  Then $$\ti \g=c^{-1}\ti E(x,t,0) a \ti E(x,t,0)^{-1} c$$ is
a new solution of the curve
flow \refep{}.  For example, if we start with the constant solution
$\g(x,t)\equiv a$, then the
corresponding solution of the $(b,j)$-flow is the vacuum solution
$u\equiv 0$.  Since the frame of $u=0$ is $E(x,t,\l)=e^{a\l x + b\l^j
t}$, the solution $\ti u$
and
$\ti\g$ are explicit given by the above formulas, which are $1$-soliton
solutions.  If we apply B\"acklund transformation repeatedly, we get the
$N$-soliton solutions
of the $(b,j)$-flow
and of the curve flow \refep{} on $M_a$.

\ms
    The algorithm of obtaining finite type solutions of the $(b,j)$-flow
is also well-known:  Fix a positive integer $k$, solve $$(\xi_1, \ldots,
\xi_k):\R^2\to
\cu_a^\perp\times \cu \times \cdots \times \cu$$
for the following two compatible ordinary differential equations:
\refeq[hw]$$\eqalign{\sum_{i=0}^k (\xi_i)_x \l^{-i}&= \left[a\l + \xi_1,\ \
\sum_{i=0}^k
\xi_i\l^{-i}\right],\cr
\sum_{i=0}^k (\xi_i)_t \l^{-i}&= \left[b\l^{j} + Q_{b,1}(\xi_1) \l^{j-1} +
\cdots + Q_{b,j}(\xi_1),\ \
\sum_{i=0}^k
\xi_i\l^{-i}\right],\cr}$$
with $\xi_0=a$.  Equate coefficients of $\l^i$ of \refhw{} to get a
system of compatible
ODE in the $x$ and $t$ variables.  So the initial value
problem for
\refhw{} has a unique
solution.
If $(\xi_1, \ldots, \xi_k)$ is a solution of \refhw{}, then
$u=\xi_1$ is a solution of the $(b,j)$-flow. Such a solution is called
a {\it finite type
solution\/}.  To obtain the corresponding solution of the curve flow
\refep{}, we first solve
$k:\R^2\to U$ for
$$k^{-1}dk = u\ dx + Q_{b,j}(u) \ dt.$$
By Proposition \refgs{} there exists a
constant $c$ such that
$\lim_{x\to -\infty} k(x,t)=c$ and
$\g= c^{-1}kak^{-1}c$ is the solution of the curve flow \refep{}
corresponding to the
finite type solution $u$ under the development map $\Phi$.

\bs

\newsection  The $U/K$-hierarchy and corresponding curve flows.\par

    Let $U$ be a compact Lie group, $\s:U\to U$ be a group involution,  $K$
the fixed point set of $\s$, and $\ck, \cp$ the $\pm 1$ eigenspaces of
$d\s_e$ on
$\cu$. Then $$\cu=\ck+\cp, \ \ [\ck,\ck]\subset \ck, \ \ [\ck,\cp]\subset
\cp, \ \ [\cp,\cp]\subset \ck.$$  In particular,
$\Ad(k)(\cp)\subset \cp$ for all $k\in K$. The quotient $U/K$ is a
symmetric space,
$\cu=\ck+\cp$ is a {\it Cartan decomposition\/} of $U/K$,  and  the
$\Ad(K)$ representation on
$\cp$ is the {\it isotropy representation\/} of $U/K$.
As reviewed in section 1,  if
$a\in \cp$, then
$\cs(\R,\cu_a^\perp\cap \ck)$ is invariant
under the odd flows in the
$U$-hierarchy, and if $k$ is odd, then the restriction of Poisson
operator $J_k$ to
$\cs(\R, \ck\cap\cu_a^\perp)$ is again a Poisson operator and the odd
flows in the
$U/K$-hierarchy are Hamiltonian with respect to $J_k$.

The curve flow corresponding to the flows in the $U/K$-hierarchy can
be obtained by
restricting the curve flow on the Adjoint orbit to a submanifold.   To
see this, let  $a\in \cp$,
$N_a$ the Ad$(K)$-orbit  at $a$ in $\cp$, and $M_a$ the Adjoint
$U$-orbit in $\cu$. Since  $N_a\subset M_a$, $C_a(\R, N_a)\subset C_a(\R,
M_a)$.

Suppose $b\in \cu_a^\perp\cap \cp$, $j$ is odd, and $u:\R^2\to
\ck\cap\cu_a^\perp$ is a solution of the
$(b,j)$-flow in the $U/K$-hierarchy.  Let $\o_\l$ be the Lax pair
\refex{} of the
$(b,j)$-flow.  Since $j$ is odd and $u\in \ck$,  $Q_{b,j}(u)\in \ck$.  Hence
$$u\ dx+ Q_{b,j}(u)\ dt$$ is a
$\ck$-valued flat connection
$1$-form.   So  $k, c$ of Proposition \refgs{} lie in $K$ and
$\g=c^{-1}kak^{-1}c$ lies in the submanifold $N_a$ of $M_a$.   Hence we have

\refclaim[hd] Theorem.  Let $U/K$ be a symmetric space, $\cu=\ck+\cp$ a Cartan
decomposition, and $a\in \cp$ such that the $\Ad(K)$-orbit $N_a$ at
$a$ in $\cp$ is a
principal orbit. Let $M_a$ be the $\Ad(U)$-orbit at $a$ in $\cu$.  Then
\item {(i)}the development map $\Phi$ maps $C_a(\R,N_a)$
isomorphically onto $\cs(\R,
\ck\cap \cu_a^\perp)$,
\item {(ii)} if $j$ is odd, then the curve flow \refep{} on
$C_a(\R,M_a)$ leaves
$C_a(\R,N_a)$ invariant, and the restriction of \refep{} to
$C_a(\R,N_a)$ corresponds to
the $(b,j)$-flow in the
$U/K$-hierarchy.  Moreover, if $k$ is odd, the restricted flow on
$C_a(\R,N_a)$ is
Hamiltonian with respect to $i^*(L_k)=i^*\Phi^*(J_{k+2})$, where
$i:C_a(\R,N_a)\to
C_a(\R,M_a)$ is the inclusion.\ei

\refpar[iq] Example.  Let $\cu=\ck+\cp$ be a Cartan decomposition of
the rank $k$
symmetric space $U/K$, $a\in \cp$ such that $N_a$ is a principal
$K$-orbit in $\cp$, and
$b\in
\cu_a$.  It is known (cf.~[Te1]) that there exists a polynomial
$f:\cp\to \R$ invariant under the
Ad$(K)$-action such that $\K f(a)=b$.  Since $f$ is $K$-invariant,
$\K f\n N_a$ is a
parallel normal field.    So
the curve flow $\g_t= (\hat b(\g))_x$ becomes
\refeq[ip]$$\g_t= (\K f(\g))_x$$ on $N_a$.
Since \refip{} corresponds to the $(b,1)$-flow in the $U/K$-hierarchy
under $\Phi$ and
the $(b,1)$-flow is completely integrable with respect to $J_3$,
\refip{} is completely
integrable with respect to the Poisson operator $\L$.  By Proposition
\refhe{}, the
restriction of $H_b$ to $C_a(\R,N_a)$ is the Hamiltonian of \refip{}
with respect to $\L$.
The higher order conserved quantities of \refip{} are the
restrictions of $F_{b,2k+1}\circ
\Phi$ to $C_a(\R,N_a)$.  Moreover, the construction of B\"acklund
transformations and finite
type solutions work in a similar way by requiring $f(\l)$ and
$\xi(\l)=\sum_{i=0}^N
\xi_i\l^{-i}$ in section 3 to satisfy the the extra reality
condition given by the involution $\s$:
$$\s(f(-\l))= f(\l), \quad d\s_e(\xi(-\l)) = \xi(\l).$$

\ms
\refpar[jb] Example.
Consider the rank one symmetric space
$$U/K=SO(n+1)/SO(n)=S^n.$$  The Cartan decomposition is
$so(n+1)= so(n) + \cp$, where $\ck=so(n)$, and
$$\cp=\left\{\pmatrix{0&-v^t\cr v&0\cr}\bigg|\  v\in \R^n\right\}.$$
Let $e_{ij}$ be the elementary matrix, and  
$$a=e_{21}-e_{12}.$$ 
The $SO(n)$-orbit in $\cp$ at $a$ under the isotropy representation of $S^n$ is the
standard unit sphere $S^{n-1}$.  Below we compute the third flow in the
$S^n$-hierarchy and the corresponding curve flow on $S^{n-1}$.  A direct
computation shows that
$\ck\cap\cu_a^\perp$ and $\cp\cap \cu_a^\perp$ are spanned by
   $$k_\a= e_{\a 2}-e_{2\a}, \ \ p_\a=e_{\a 1} - e_{1 \a}, \ \ 3\leq \a\leq n+1$$
respectively.  
Henceforth in this example, we assume $3\leq \a,\b \leq n+1$.     Let
$u=\sum_{\a=3}^{n+1} u_\a k_\a$, and $Q_j= Q_{a,j}(u)$. We use \refad{} and
\refiw{} to compute the $Q_j$'s.  Note that
$Q_0=a$ and $Q_1=u$.
The reality conditions for the $S^n$-hierarchy are
\refeq[jc]$$\overline{A(\bar\l)} = A(\l), \quad A(\l)^t + A(\l)=0,
\quad I_{n,1}^{-1}
A(-\l) I_{n,1} = A(\l),$$ where $I_{n,1}=\diag(1, \ldots, 1, -1)$.
It is known (cf.~[TU1]) that
$\sum_{j=0}^\infty Q_j
\l^{-j}$ satisfies the reality conditions \refjc{}.  Hence
   $Q_i$ is in $\ck$ for odd $i$ and in $\cp$ for
even $i$.  In particular, $Q_2\in \cp$.  Write $Q_2= y_0 a +
\sum_{\a=3}^{n+1}  y_\a p_\a$.  Then  
$$[Q_2,a] = (Q_1)_x+[u,Q_1]= u_x$$
implies that
$$y_\a= -(u_\a)_x.$$ To compute $y_0$, we use condition
\refiw{} to conclude
that $$\tr(a+u\l^{-1}+ Q_2\l^{-2} + \cdots)^2 \sim \tr(a^2)$$ as an asymptotic
expansion.  Compare the coefficients of
$\l^{-2}$ in the above equation to get
$$\tr(aQ_2+Q_2a+u^2)=0.$$
This implies $y_0= -{\N u\N^2\over 2}$.  Hence
$$Q_2= -{\N u\N^2\over 2}\ a - \sum (u_\a)_x p_\a.$$

\ms
We know $Q_3\in \ck$, and
\refeq[ad]$$(Q_2)_x+[u,Q_2]=[Q_3,a].$$  Let $Q_3=\sum_{i,j=1}^n y_{ij}
e_{ij}$. Then \refad{} implies that
$$y_{2\a}= (u_\a)_{xx} +{\N u\N^2\over 2} u_\a.$$
Compare the coefficients of $e_{\a\b}$ for
$\a,\b\ge 3$ in
$$(Q_3)_x+[u,Q_3]=[Q_4,a]$$ to get
$$(y_{\a\b})_x= -u_\a(u_\b)_{xx} + u_\b(u_\a)_{xx}.$$
The right hand side is equal to $((-u_\a(u_\b)_x +
u_\b(u_\a)_x)_x$. By \refjd{}, $Q_{b,j}(u)\in \cs(\R, \cu)$.  So
$$y_{\a\b}= -u_\a(u_\b)_x + u_\b(u_\a)_x.$$
Hence
\refeq[iy]$$Q_3=\sum_\a \left(-(u_\a)_{xx} -{\N u\N^2\over 2}
u_\a\right)k_\a + \sum_{\a,\b}
y_{\a\b} e_{\a\b}.$$
Now we compute the third flow
\refeq[hz]$$u_t=(Q_3)_x+[u,Q_3]$$ directly. The coefficient of
$e_{2\a}$ of the left
hand side of
\refhz{} is $-(u_\a)_t$, and of the right hand side is
$$\eqalign{& (y_{2\a})_x+ \sum_\b u_{2\b} y_{\b\a} -\sum_\b y_{2\b}
u_{\b\a}\cr &=  \left((u_\a)_{xx} +{\N u\N^2\over 2} u_\a\right)_x  - \sum_\b u_\b
(-u_\b(u_\a)_x + u_\a(u_\b)_x)\cr
&= (u_\a)_{xxx} +{\N u\N^2\over 2} (u_\a)_x +\sum_\b u_\b(u_\b)_xu_\a
+ \N u\N^2 (u_\a)_x
- u_\a \sum_\b u_\b (u_\b)_x .\cr}$$
So the third flow is the vector modified KdV equation:
\refeq[]$$(u_\a)_t=  -(u_\a)_{xxx} -{3\over 2} \N u\N^2 (u_\a)_x,$$
or equivalently
$$u_t= -\left(u_{xxx}+ {3\over 2} \N u\N^2 u_x\right), \eqno({\rm vmKdV\/})$$
where $u:\R^2\to \R^{n-1}$.

The Ad$(K)$-orbit at $a$ in $\cp$ is the standard sphere $S^{n-1}$ in
$\cp$.  By
Proposition \refhs{},  the curve flow on $S^{n-1}$ corresponding to
the first flow is
the translation flow $\g_t= \g_x$.  The curve flow on $S^{n-1}$
corresponding to the
third flow in the
$S^n$-hierarchy is
\refeq[ix]$$\g_t= g[Q_3(u),a]g^{-1},$$ where $\g= gag^{-1}$ and $g^{-1}g_x= u$.
Substitute \refiy{} into the right hand side of \refix{} to get
\refeq[iz]$$\g_t= -\sum_{\a=3}^{n+1} \left((u_\a)_{xx} + {\N u\N^2\over
2} u_\a\right) \
gp_\a g^{-1}.$$
A direct computation gives
$$\eqalign{\g_x&=  \sum_{\a}u_\a gp_\a g^{-1},\cr
  \g_{xx} &= \sum_{\a}(u_\a)_x gp_\a g^{-1} - \Vert u\Vert^2 gag^{-1},\cr
  \g_{xxx}
  &=\sum_\a((u_\a)_{xx} - \Vert u\Vert^2 u_\a) gp_\a g^{-1} - 3\sum_\a u_\a
(u_\a)_x
  gag^{-1}.\cr}$$
We can express the right hand side of \refiz{} in terms of $\g$ and its
$x$-derivatives to get
$$\g_t = -\left(\g_{xxx} + 3\li \g_x, \g_{xx}\ri \ \g + {3\over 2} \N
\g_x\N^2 \g_x\right).$$

\bs

\newsection  Weakly non-linear hydrodynamic systems.\par

Let $M$ be a principal orbit  of the isotropy
representation of a symmetric space $U/K$, and $v$ a parallel normal field.
Ferapontov noted in [F2] and [F4] that the curve flow \refel{}
$$\g_t= (v(\g))_x = -A_{v(\g)}(\g_x)$$
  is a  weakly non-linear hydrodynamic system on
$M$. In this section, we study this curve flow on
principal orbits of the isotropy representation of a symmetric space
as a hydrodynamic
system.

First we review some definitions and results on hydrodynamic
systems. Given a
smooth map
$v=(v_{ij}):\R^n\to gl(n,\R)$, the first order quasilinear system
for $u=(u_1, \ldots, u_n)^t:\R^2\to \R^n$,
\refeq[bs]$$(u_i)_t= \sum_{j=1}^n v_{ij}(u)(u_j)_x,$$
is called a {\it hydrodynamic system\/} on $\R^n$.
The system \refbs{} is said to be {\it diagonalizable\/}
in an open subset $\co$ of $\R^n$ if there exist local coordinates
$(y_1, \ldots, y_n)$ in
$\co$ such that the system
\refbs{} is of the form
$$(y_i)_t= v_i(y) (y_i)_x.$$ It is called {\it weakly non-linear\/} or {\it
linearly degenerate\/} on $\co$ if
\item {(i)} the eigenvalues of $(v_{ij}(x))$ have
constant  multiplicities $m_1, \ldots, m_s$
(so the corresponding eigenvalue functions $\l_1, \ldots, \l_s$ are smooth),
\item {(ii)}  if $\xi$ is an eigenvector of $v=(v_{ij})$ with 
eigenvalue $\l_j$,
then
$d\l_j(\xi)=0$.

\ms
A {\it hydrodynamic system\/} on a manifold $M$ is
a first order quasilinear system for $\g:R^2\to M$ of the form:
\refeq[as]$$\g_t= P(\g)(\g_x),$$ where $P$ is a smooth section of $L(TM,TM)$.
    It is easy to see that this system in local coordinates looks like a
hydrodynamic
system on an open subset of $\R^n$.   System \refas{} is called {\it
diagonalizable\/} ({\it weakly non-linear\/} respectively) if locally
it is diagonalizable (weakly
non-linear  respectively).

Recall that
a submanifold $M$ of $\R^{n+k}$ is called {\it isoparametric\/} if the
normal bundle is flat
and the principal curvatures along any parallel normal field are constant. By
definition of weak non-linearity,
the flow
\refel{} on an  isoparametric submanifold is a
weakly non-linear hydrodynamic system.  Since a principal orbit $M$ of the
isotropy representations of a symmetric space is isoparametric,
\refel{} on $M$ is a weakly
non-linear hydrodynamic system.

     Let $\R^n$ be equipped with the standard inner product, and $\K$
the Levi-Civita
connection.  Then $\K_{u_x}$ is a Poisson operator on $\cs(\R,\R^n)$.
     Given a smooth function $f:\R^n\to \R$, the Hamiltonian equation of
the functional
$$F(u)=\int_{-\infty}^\infty f(u(x))dx$$ with respect to $\K_{u_x}$  is
\refeq[bu]$$u_t= \K_{u_x} (\K f(u)).$$
Write \refbu{} in the standard coordinates of $\R^n$ to get
$$u_t= (\K f(u))_x, $$ i.e.,
\refeq[is]$$(u_i)_t= \sum_{j=1}^n f_{u_i u_j}(u)(u_j)_x, \quad 1\leq
i\leq n. $$
So \refis{} is a hydrodynamic Hamiltonian system on $\R^n$.   Dubrovin and
Novikov investigated the systems \refis{} on $\R^n$,
and obtained many  remarkable results (cf.~[DN]).  Novikov conjectured that if
\refbu{} is diagonalizable then it is  a completely integrable
Hamiltonian system. Tsarev
proved this conjecture and gave a complete classification of such
systems (see [Ts]).
Moreover, all constant of motions of \refbu{} are of zero order.  We
remark that the boundary
conditions of the Poisson operator $\K_{u_x}$ are not taken into
account in these results.

Let $M$ be a principal orbit in $\cp$ of the isotropy representation
of the symmetric space
$U/K$.  We have seen in Example \refiq{} that the curve flow \refel{}
is of the form
\refeq[hv]$$\g_t= (\K f(\g))_x$$
on $M$ for some $K$-invariant polynomial $f:\cp\to \R$.  Note that
\refhv{} on $\cp$ is a
hydrodynamic system and is Hamiltonian with respect to
$\K_{u_x}$.  Although $C_a(\R, M)$ is invariant under the flow
\refhv{}, the restriction of
$\K_{u_x}$ to $C_a(\R,M)$ is not a Poisson structure.  However, we
have shown that
\refhv{} is a completely integrable Hamiltonian system with respect
to the Poisson operator
$\L$ on $C_a(\R,M_a)$, its Hamiltonian is of zero order (given by the
height function
$H_b$), and it has infinitely many higher order conserved quantities.
Ferapontov ([F3]) noted that
\refhv{} is non-diagonalizable if $U/K$ is $SU(3)/SO(3)$.  Below we use
submanifold geometry to prove directly that
\refhv{} is non-diagonalizable on any irreducible isoparametric
submanifold.

First we need to review
some isoparametric theory (cf.~[Te1]).  Let $M^n\subset \R^{n+k}$ be an
isoparametric
submanifold. Then there exist smooth subbundles
$E_1, \ldots, E_p$ of $TM$ and parallel normal fields $v_1, \ldots,
v_p$ such that
\item {(i)} $TM=\oplus_{i=1}^p E_i$,
\item {(ii)} if $v$ is a parallel normal vector field on $M$ then the
shape operator
$$A_v\vert E_i=\li v, v_i\ri \id_{E_i},$$
\item {(iii)} there exists a parallel normal field $v$ such that
$$\li v, v_1\ri, \li v, v_2\ri,
\ldots, \li v, v_p\ri$$ are distinct; in particular, $E_1, \ldots,
E_p$ are eigenspaces of
$A_v$.

\ss
\ni
The $E_i$'s and $v_i$'s are called {\it curvature distributions\/}
and {\it curvature
normals\/} of
$M$.   An isoparametric submanifold $M$ of $\R^{n+k}$ is {\it
irreducible\/} if $M$ is not a
product of two lower dimensional isoparametric submanifolds.  It is known that
\item {(i)} if $\li v_i, v_j\ri=0$ for all $i\not=j$,
then $M$ is a product of spheres,
\item {(ii)} principle orbits of isotropy representations of
irreducible symmetric spaces are
irreducible isoparametric submanifolds.

\refclaim[be] Proposition.  Let $X:M^n\to \R^{n+k}$ be a compact,
irreducible, isoparametric
submanifold with $k\geq 2$,   $v_1, \ldots, v_p$ the curvature
normals of $M$, and $v$ a
parallel normal field on $M$ such that
$\li v_1, v\ri, \ldots, \li v_p,v\ri$  are distinct. Then the
hydrodynamic system
$\g_t=(v(\g))_x$ on $M$  is not diagonalizable.

\proof  Suppose $\g_t= (v(\g))_x=-A_{v(\g)}(\g_x)$ is diagonalizable.
Then there is a local
coordinate system $y$ such that
$$\eqalign{&I = \sum_{i=1}^n a_i(y)^2 dy_i^2,\cr
&A_v(X_{y_i}) =  \l_i  X_{y_i}, \qquad 1\leq i\leq n\cr}$$ for some smooth
functions $a_i$ and constants $\l_i$.
Set $e_i= X_{y_i}/a_i$. Then $e_1, \ldots, e_n$ is a local orthonormal
tangent frame on $M$.  Since $E_1,
\ldots, E_p$ are eigenspaces of $A_v$, we may arrange the indices $i$ so that
$$\cases{e_1, \ldots, e_{i_1}& span $E_1$\cr
e_{i_1+1}, \ldots, e_{i_2} & span $E_2$, \cr
    \ldots, & $\ldots$\cr
e_{i_{p-1}+1}, \ldots, e_n & span $E_p$.\cr}$$
Set $$I_j=\{i_{j-1}+1, \ldots, i_j\}.$$
The dual coframe of $e_i$ is
$w_i=a_i \ dy_i$. Let $e_{n+1}, \ldots, e_{n+k}$ be a parallel normal frame on
$M$.  Write
$$de_i =  \sum_{j=1}^{n+k} w_{ji}e_j, \quad 1\leq i\leq n+k.$$
It follows from elementary local submanifold geometry that we have
$$\eqalign{w_{\a\b}& =0,  \quad n+1\leq \a, \b\leq n+k\cr
    w_{i\a}&=\l_{i\a} w_i, \quad 1\leq i\leq n,\ \  n+1\leq \a\leq n+k, \cr
w_{ij}&= {(a_i)_{y_j}\over a_j }\ dy_i -  {(a_j)_{y_i}\over a_i}\
dy_j, \quad 1\leq i\not=j
\leq n,\cr
dw_{AB}&= -\sum_{C=1}^{n+k} w_{AC}\wedge w_{CB}, \quad 1\leq A,B\leq n+k,\cr}$$
and
$$v_m=\sum_{\a=n+1}^{n+k} \l_{i\a} e_\a, \ \ {\rm where\ } i\in I_m$$
The  Codazzi equation
$dw_{i\a}=\sum_{j=1}^n w_{ij}\wedge w_{j\a}$ implies that
\refeq[hp]$${(a_i)_{y_j}\over a_j} = {(\l_{i\a} a_i)_{y_j}\over
\l_{j\a} a_j}= {\l_{i\a}\over
\l_{j\a}} \ {(a_i)_{y_j}\over a_j}$$  for all
$1\leq i\not=j\leq n$ and $n+1\leq \a\leq n+k$.  If
$e_i\in E_k$, $e_j\in E_m$ and $k\not=m$, then since $v_1, \ldots,
v_p$ are distinct, there
exists $\a$ such that $\l_{i\a}\not= \l_{j\a}$.  It follows from
equation \refhp{} that
\refeq[ja]$${\rm If\ } r\not=s,\ i\in I_r, \ {\rm and\ } j\in
I_s{\rm, \  then\ } w_{ij}=0.$$

Next we claim that $\li v_r, v_s\ri=0$ if $ r\not=s$.  To see this,
let $i\in I_r$ and
$j\in I_s$.  The Gauss equation implies that
$$dw_{ij}=-\sum_{m=1}^n w_{im}\wedge w_{mj} + \sum_{\a=n+1}^{n+k}
w_{i\a}\wedge w_{j\a}.$$
  By \refja{}, the left hand side
and  the first term of the right hand side of the above equation are zero.   So
$$\sum_{\a=n+1}^{n+k} w_{i\a}\wedge w_{j\a} =0.$$
But $$ \sum_{\a=n+1}^{n+k} w_{i\a}\wedge w_{j\a} = \sum_\a \l_{i\a}\l_{j\a}\
w_i\wedge w_j = \li v_r, v_s\ri \ w_i\wedge w_j,$$
which proves the claim.  So
$M$ is the product  of standard spheres. This
contradicts the assumption that $M$ is irreducible.
Hence
$\g_t=(v(\g))_x$ is
not diagonalizable. \qed

\vfil\eject
 \hsize  5.875 true  in 
\hoffset=  42 true pt

      \vsize  9.1 true  in

    \Bibliography

\a //DN//Dubrovin, B.A., Novikov, S.P.//Hamiltonian formalism of
one-dimensional systems of hydrodynamic type, and the Bogolyubov-Whitham
averaging method//Soviet Math. Dokl.//27//1983//781-785////

\b //FT//Faddeev, L.D., Takhtajan, L.A.//Hamiltonian methods in the theory
of solitons//Springer-Verlag////1987//////

\a //F1//Ferapontov, E.V.//Dirac reduction of the Hamiltonian
operator $\d^{ij}
{d\over dx}$ to a submanifold of the Euclidean space with flat normal
bundle//Funk. Anal. Pril.//26//1992//83-86////

\a //F2//Ferapontov, E.V.//On the matrix Hopf equations and
integrable hamiltonian
systems of hydrodynamic type, which do not posses Riemann
invariants//Physics Lett. A//179//1993//391-397////

\a //F3//Ferapontov, E.V.//On integrability of $3\times 3$ semihamiltonian
hydrodynamic type systems which do not  posses Riemann
invariants//Physics D//63//1993//50-70////

\a //F4//Ferapontov, E.V.//Isoparametric hypersurfaces in spheres, integrable
non-diagonalizable systems of hydrodynamic type, and $N$-wave
systems//Diff. Geom. Appl.//5//1995//335-369////

\b //H//Helgason, S.//Differential Geometry, Lie groups, and symmetric
spaces//Academic press, NY////1978//////

\a //Sa//Sattinger, D.H.//Hamiltonian hierarchies on semi-simple Lie
algebras//Stud. Appl. Math.//72//1984//65-86////

\a //Te1//Terng, C.L.//Isoparametric
submanifolds and their Coxeter groups//J.~Differential
Geometry//21//1985//79--107////

\a//Te2//Terng, C.L.//Solitons and Differential Geometry//J. of
Differential Geometry //45//1997//407-445////

\a //TU1//Terng, C.L., Uhlenbeck, K.//B\"acklund transformations
and loop group actions//Comm. Pure. Appl. Math.//53//2000//1-75////

\p //TU2//Terng, C.L., Uhlenbeck, K.//Schr\"odinger flow on
Grassmannians//////////math.DG/9901086, to appear in ``Integrable systems,
geometry and topology'' published by International Press.//

\a //Ts//Tsar\"ev, S.P.//The geometry of Hamiltonian systems of
hydrodynamic type. The generalized Hodograph
method//Math. USSR Izvestiya//37//1991//397-419////

\b //W//Wolf, J.A.//Spaces of Constant Curvature,
3rd Edition//Publish or Perish, Boston////1972//////

\bs\bs
\ni Chuu-Lian Terng: Department of Mathematics, Northeastern
University, Boston, MA  02115.

\ni Email: terng@neu.edu
\ss
\ni Gudlaugur Thorbergsson: Mathematisches Institut der
Universit\"at zu K\"oln, Weyertal
86--90, D-50931 K\"oln, Germany. 
Email: gthorbergsson@mi.uni-koeln.de
\bs


\end